\documentclass[11pt,reqno,UTF8,twoside]{amsart}
\usepackage{mathrsfs}
\usepackage{amsfonts,amssymb,amsmath,amsthm}
\usepackage{cite}
\usepackage[colorlinks=true,citecolor=red]{hyperref}
\usepackage{titletoc}
\usepackage{geometry}
\usepackage{bm}
\usepackage{indentfirst}
\usepackage{graphicx}
\usepackage{float} 
\usepackage{booktabs}
\usepackage{longtable}
\usepackage{subfigure}
\usepackage{stmaryrd}
\usepackage{enumerate}
\usepackage{tikz}
\usepackage{ulem}
\usepackage{color,soul}
\usepackage{setspace}
\usepackage{stmaryrd}
\usepackage[pagewise]{lineno} 
\allowdisplaybreaks[2]
\usepackage{appendix}
\usepackage{cases}
\usepackage{todonotes}
\usepackage{enumitem}
\usepackage{comment}
\usepackage{graphicx}

\numberwithin{equation}{section}
\newtheorem{theorem}{Theorem}[section]
\newtheorem{lemma}[theorem]{Lemma}
\newtheorem{corollary}[theorem]{Corollary}
\newtheorem{proposition}[theorem]{Proposition}

\newtheorem{definition}[theorem]{Definition}
\theoremstyle{remark}

\newtheorem{rmk}[theorem]{Remark}

  \href{}{}

\hypersetup{linkcolor=blue}

\makeatletter
\@namedef{subjclassname@2020}{\textup{2020} Mathematics Subject
Classification}
\makeatother

\begin{document}
    
\title[Mixing for vortex model]{Exponential Mixing for the Randomized Two-vortex Model}

%\author[Y. Liu]{Yong Liu}
%\address[Yong Liu]{LMAM and School of Mathematical Sciences, Peking University, 100871, Beijing, China.}
%\email{liuyong@math.pku.edu.cn}

\author[Y. Shi]{Yankai Shi}
\address[Yankai Shi]{School of Mathematical Sciences, Peking University, 100871, Beijing, China.}
\email{shiyankai@stu.pku.edu.cn}

%\author[B. Tang]{Bin Tang}

	

\begin{abstract}
This paper investigates the mixing of passive scalars advected by a randomized two-vortex model on $\mathbb{T}^2$, in which the vortices interact through a regularized Biot–Savart kernel and are driven by independent Brownian motions. We establish quenched exponential mixing for the Lagrangian flow, characterized by the exponential correlation decay of mean-zero $H^s$ observables. The proof proceeds by establishing ergodicity of the associated two-point Markov process. We combine a uniform contraction estimate with controllability arguments to verify the Lyapunov drift and minorization conditions, allowing us to apply Harris-type theorem.

\end{abstract}

\maketitle
\setcounter{tocdepth}{1}
\tableofcontents

%\tableofcontents
\section{Introduction}\label{sec: Introduction}
Passive scalar mixing by incompressible flows is of fundamental importance in fluid mechanics and the natural sciences. Physically, through repeated stretching and folding, advection transfers scalar fluctuations from macroscopic spatial scales to microscopic ones. Subsequently, diffusion acts effectively on these fine scales, ultimately driving the system toward equilibrium. Therefore, the study of mixing is essential for understanding turbulent and chaotic dynamical systems.

\vspace{0.5em}

In this paper, we investigate the transport equation
\begin{equation}\label{eq: advection-diffusion}
    \begin{cases}
        \partial_t\rho+u\cdot\nabla\rho=0,\\
        \rho\big|_{t=0}=\rho_0,
    \end{cases}
\end{equation}posed on the two-dimensional torus $\mathbb{T}^2=\mathbb{R}^2/\mathbb{Z}^2$. Here, $u\colon [0,+\infty)\times\mathbb{T}^2\rightarrow\mathbb{R}^2$ is a divergence-free velocity field. The solution $\rho$ to \eqref{eq: advection-diffusion}, representing the concentration of a passive scalar, is assumed to be mean-zero:
\begin{equation*}
    \int_{\mathbb{T}^2}\rho(t,x)\mathrm{d}x=0,\qquad \forall t\geq 0.
\end{equation*}As proposed in \cite{MMP-05,Thiffeault-11}, a conventional measure to quantify the mixing of $\rho$ is the homogeneous negative Sobolev norms
\begin{equation*}
    \|\rho(t)\|_{\dot{H}^{-s}}^2:=\sum_{k\in\mathbb{Z}^2\setminus\{\mathbf{0}\}}\left|k\right|^{-2s}\left|\hat{\rho}_k(t)\right|^2, \qquad s>0,
\end{equation*}whose decay indicates that the $L^2$-energy of $\rho$ cascades to high-frequency modes.

Alternatively, \eqref{eq: advection-diffusion} can be reformulated within the Lagrangian framework generated by $u$:
\begin{equation*}
    \begin{cases}
       \mathrm{d}{\phi}_t(x) = u(t,\phi_t(x))\mathrm{d}t,\\
        \phi_0(x)=x.
    \end{cases}
\end{equation*}Here, $\phi_t(x)$ traces the trajectory of a passive scalar particle initialized at $x$ and satisfies 
\begin{equation*}
    \rho(t,x)=\rho_0\circ\phi_t^{-1}(x).
\end{equation*}We adopt this Lagrangian perspective hereafter to study the mixing properties of a randomized two-vortex model, introduced below.

\subsection{Model and main result}\label{subsec: Model and main result}
Several recent works have constructed exponential mixers based on explicit steady solutions to the Euler equations, including randomized alternating shear flows \cite{BCZG-23}, cellular flows with randomly moving centers \cite{NFS-26}, and randomized ABC flows \cite{CZNF-26}. Time dependence in these models is introduced through random switching, translations, or parameter changes within a prescribed velocity fields. A broader question therefore arises: whether analogous mixing estimates hold for more general Euler-type dynamics under suitable random perturbations, without an explicitly designed stirring protocol?

The point-vortex system provides a natural finite-dimensional model for investigating this question. As a classical approximation of Euler flows \cite{MP-93,Aref-07}, the system evolves according to
\begin{equation*}
    \mathrm{d}X_i=\frac{1}{N}\sum_{j\neq i}K(X_i-X_j)\mathrm{d}t,\qquad i=1,2,\cdots,N,
\end{equation*}where the Biot–Savart kernel is
\begin{equation*}
    K(x):=\sum_{k\in\mathbb{Z}^2\setminus\{\mathbf{0}\}}\frac{ik^\perp}{2\pi|k|^2}e^{2\pi ik\cdot x}.
\end{equation*}Unlike prescribed stirring mechanisms, each vortex is advected by the mutual interactions within the system. Chaotic advection generated by deterministic point-vortex motions has been studied extensively; see \cite{Aref-84,BSA-03}. Nevertheless, numerical observations in \cite{BBPV-94} revealed that regular islands may persist around vortex centers, where the passive particle dynamics is characterized by a null Lyapunov exponent. This coexistence indicates that chaotic vortex motion does not necessarily ensure global mixing.

\vspace{0.5em}
Motivated by the limitations of the deterministic regime, we introduce a randomized point-vortex model. To regularize the singularity of $K$, we fix $\epsilon>0$ and employ the following Gaussian-regularized periodic kernel:
\begin{equation}\label{eq: regularized periodic kernel}
    K_\epsilon(x):=\sum_{k\in\mathbb{Z}^2\setminus\{\mathbf{0}\}}\frac{e^{-4\pi^2\varepsilon^2 |k|^2}}{2\pi|k|^2}ik^\perp e^{2\pi ik\cdot x},\qquad k^\perp=(-k_2,k_1),
\end{equation}The model considered in this paper is the randomly perturbed two-vortex system
\begin{equation}\label{eq: vortices}
    \begin{cases}
    \mathrm{d}X_1(t)=K_\epsilon\!\left(X_1(t)-X_2(t)\right)\mathrm{d}t+\mathrm{d}B_1(t),\quad X_1(0)=X_1,\\
    \mathrm{d}X_2(t)=K_\epsilon\!\left(X_2(t)-X_1(t)\right)\mathrm{d}t+\mathrm{d}B_2(t),\quad X_2(0)=X_2,
    \end{cases}
\end{equation}
where $B_1,B_2$ are independent Brownian motions on $\mathbb{T}^2$. 

Let $Z=(X_1,X_2)\in\mathbb{T}^4$ denote the joint vortex state, and define
\begin{equation*}
    b(Z):=\left(K_\epsilon\!\left(X_1-X_2\right),K_\epsilon\!\left(X_2-X_1\right)\right).
\end{equation*}The vortex dynamics can be rewritten as 
\begin{equation*}
    \mathrm{d}Z_t=b(Z_t)\mathrm{d}t+\mathrm{d}B_t,
\end{equation*}which induces the corresponding random velocity field
\begin{equation}\label{eq: velocity field}
    u(t,x)=u\!\left(Z_t,x\right):=K_\epsilon\!\left(x-X_1(t)\right)+K_\epsilon\!\left(x-X_2(t)\right).
\end{equation}Since $K_\varepsilon$ is smooth and divergence-free, for any initial configuration $Z$, this velocity field $u$ generates a volume-preserving Lagrangian flow
\begin{equation}
    \mathrm{d}\phi_t^Z(x)=u\!\left(Z_t,\phi_t^Z(x)\right)\mathrm{d}t.\label{eq:Lagrangian flow}
\end{equation}With these preparations, we are ready to state our main result.
\begin{theorem}\label{thm: exponential mixing}
    Let $\phi_t$ denote the Lagrangian flow of equation \eqref{eq:Lagrangian flow}. For any $q,s>0$, there exist a constant $\gamma_{s,q}>0$ and a random variable $D=D\!\left(\omega,Z\right)\colon\Omega\times\mathbb{T}^4\rightarrow [0,+\infty)$ such that for all mean-zero observables $f,g\in H^s\!\left(\mathbb{T}^2\right)$, the following holds
    \begin{equation*}
        \left|\int_{\mathbb{T}^2}f(x)g\!\left(\phi_t(x)\right)\mathrm{d}x\right|\leq D\|f\|_{H^{s}}\|g\|_{H^s}e^{-\gamma_{s,q}t},
    \end{equation*}almost surely for all $t>0$. Furthermore, there exists a deterministic constant $D_{s,q}$ such that
    \begin{equation*}
        \mathbb{E}[D^q]\leq D_{s,q}.
    \end{equation*}
\end{theorem}
\begin{rmk}\label{rmk: N vortex}
    The approach developed in the proof of Theorem \ref{thm: exponential mixing} can be directly extended to randomized $N$-vortex systems for $N \geq 2$, subject to a necessary modification of the exponential mixing rate $\gamma_{s,q}$.
\end{rmk}

Theorem \ref{thm: exponential mixing} is a quenched correlation decay estimate. By volume preservation and duality, this indicates the exponential mixing for the transport equation \eqref{eq: advection-diffusion} with $u$ given by \eqref{eq: velocity field}. We emphasize that all constants herein may depend on the regularization parameter $\epsilon$.

Our long-term objective is to understand mixing in the unregularized point-vortex system as $\epsilon\rightarrow 0$, and eventually in the dynamical limit $N\to\infty$. The present work provides a first step toward these questions, with the regularization scale $\epsilon>0$ fixed.

A key challenge is to quantify the dependence of the mixing rate on $\epsilon$ and $N$. Within the Harris framework (see Theorem \ref{thm: Discrete Harris theorem}), this calls for more precise estimates of the Lyapunov drift and minorization constants, particularly tracking the accessibility times and transition probability lower bounds. Such quantitative estimates provide a basis for understanding whether exponential mixing persists in these limiting regimes.

\subsection{Proof Strategy} The proof of Theorem \ref{thm: exponential mixing} builds on the methodology developed in \cite{DKK-04,BBPS-mixing,BCZG-23,LS-26}. Following \cite{LS-26}, we construct an explicit Lyapunov function via a uniform contraction estimate, without resorting to spectral perturbation of the twisted projective semigroup. We now briefly outline the main ingredients of the argument.

By Sobolev duality of $\|\cdot\|_{H^{-s}}$ and the identity $\rho_t=\rho_0\circ\phi_t^{-1}$, we obtain
\begin{equation*}
    \|\rho_t\|_{H^{-s}}=\sup_{\|\varphi\|_{H^s}=1}\int_{\mathbb{T}^2}\varphi(x)\rho_t(x)\mathrm{d}x=\sup_{\|\varphi\|_{H^s}=1}\int_{\mathbb{T}^2}\rho_0(x)\varphi\!\left(\phi_t(x)\right)\mathrm{d}x.
\end{equation*}Therefore, the desired exponential mixing can be formulated in terms of quenched correlation decay. Following the standard argument (see Section \ref{subsec: Proof of Theorem 1.1}), it suffices to establish exponential ergodicity of the two-point Markov process
\begin{equation*}
    (Z_t,\phi_t(x),\phi_t(y)),\qquad x\neq y.
\end{equation*}Since its state space excludes the diagonal and is therefore non-compact, we apply Harris theorem, reducing the ergodicity proof into two ingredients: a Lyapunov drift condition and a minorization condition.

\subsubsection{Lyapunov drift condition} The key ingredient is the uniform contraction condition employed in \cite{LS-26}:
\begin{equation*}
    \sup_{(Z,x,v)\in\mathbb{T}^4\times\mathbb{T}^2\times S^1}
    \mathbb{E}_{Z,x}\|D_x\phi_T(x)v\|^{-p}
    \leq\gamma',
\end{equation*}
for some $T>0$ and $p,\gamma'\in(0,\frac12)$.

To this end, we first prove the positivity of the top Lyapunov exponent $\lambda_1$ via the Furstenberg’s criterion formulated in \cite{BBPS-Lyapunov}. The Furstenberg--Khasminskii formula, combined with uniform Ces\`aro convergence then yields uniformly positive expected logarithmic stretching, establishing the desired contraction. 

Following \cite{LS-26}, we consider the explicit Lyapunov function
\begin{equation*}
    V(Z,x,y)=\operatorname{dist}_{\mathbb{T}^2}(x,y)^{-p},\qquad x\neq y.
\end{equation*}Uniform derivative bounds allow us to compare nearby trajectories with the linearized flow, yielding contraction of $V$ near the diagonal. Away from the diagonal, the inverse Lipschitz estimate provides a uniform bound. Combining these estimates gives the global drift condition
\begin{equation*}
    P_T^{(2)}V\leq\gamma V+C, \qquad \gamma\in(0,1).
\end{equation*}

\subsubsection{Minorization condition} Although the noise acts only on
the vortex variables, the spanning property in Lemma \ref{lem: nondegeneracy lemma} implies the parabolic H\"ormander condition for the two-point process. This yields smooth transition densities and the strong Feller property.

We then establish topological irreducibility through global approximate controllability; see, e.g., \cite{Jurdjevic-97}. The key observation is that for any $x\neq y$, the convex hull of the attainable velocities contains a neighborhood of the origin. This allows us to design admissible control paths $h_1,h_2$ that effectively steer the two-point system arbitrarily close to any prescribed target configuration. Combining irreducibility with smooth density estimates, we obtain the minorization condition on each prescribed compact sublevel set of $V$, at a skeleton time $m\tau$.

\subsection{Discussion} The mathematical theory of mixing in deterministic systems is well developed. For incompressible velocity fields $u \in L_t^\infty W_x^{1,p}$ with $p>1$, mixing cannot proceed faster than exponentially \cite{CDL-08,IKX-14,Seis-13}, while the critical case $p=1$ remains open. Exponentially mixing flows $u$ have been constructed for prescribed initial data \cite{ACM-19A,YZ-17} and as universal mixers for $\rho_0\in L^\infty\cap H^\sigma$; see \cite{EZ-19}. 

In the stochastic setting, foundational results \cite{BS-88,DKK-04} demonstrated almost-sure exponential mixing for incompressible SDEs. In recent years, Bedrossian et al. \cite{BBPS-mixing} proved quenched exponential mixing with $u$ advected by stochastic Navier--Stokes equations. Their proof draws on Lagrangian chaos \cite{BBPS-Lyapunov} and the ergodic theory of SPDEs \cite{HM-06,HM-11EJP}. This result was later extended to finite-dimensional degenerate forcing \cite{CR-26}. Mixing has also been actively studied in random dynamical systems, including the randomized Pierrehumbert model \cite{BCZG-23,CIS-25,son-25}, ABC flows \cite{CZNF-26} and the Chirikov standard map \cite{LS-26}.

For point-vortex models, deterministic results \cite{BBPV-94} have investigated chaotic advection and the persistence of regular regions. The literature on stochastic point-vortex systems has largely addressed well-posedness and hydrodynamic limits. Flandoli et al. \cite{FGP-11} proved global well-posedness with nondegenerate transport noise. In the limiting regime, propagation of chaos toward the two-dimensional Navier--Stokes equations was established in \cite{FHM-14}, while \cite{JW-18} obtained quantitative estimates for interacting stochastic particles with singular kernels. Although these works bridge stochastic vortex dynamics and continuum fluid equations, Theorem \ref{thm: exponential mixing} provides, to our knowledge, the first rigorous proof of exponential mixing driven by a vortex model.

\subsubsection{Extensions and outlook} The present work opens several avenues for future investigation. First, it is natural to consider the advection-diffusion equation
\begin{equation*}
    \begin{cases}
        \partial_t\rho+u\cdot\nabla\rho=\kappa\Delta\rho,\\
        \rho\big|_{t=0}=\rho_0,
    \end{cases}
\end{equation*}and study the associated enhanced dissipation. Following \cite{BBPS-enhanced}, almost-sure enhanced dissipation reduces to the $\kappa$-uniform ergodicity of the associated two-point semigroups $P_t^{(2),\kappa}$. For fixed $\epsilon>0$, the uniform minorization condition is a direct consequence of the present hypoellipticity and controllability arguments. Hence, a complete extension would require compatible uniform Lyapunov drift estimates.

A second direction concerns the singular limit $\epsilon\rightarrow0$. The principal obstruction to passing directly to the limit is the loss of the global derivative estimates used in the near-diagonal analysis. The controllability and minorization arguments must also be adapted to the unregularized dynamics. Therefore, progress toward the singular system requires a replacement for the global smooth estimates on which the present proof relies.

Finally, the fixed-$N$ extension in Remark \ref{rmk: N vortex} raises the question of how the mixing estimates behave as $N\to\infty$. This problem must be formulated under specific scalings of the circulations and noise, since the limiting fluid equation depends on both the noise structure and the chosen scaling; see \cite{FHM-14,CM-20,JW-18}. Identifying which scalings preserve mixing, and how they interact with the singular limit $\epsilon\rightarrow0$, remains a central problem.

\subsection*{Organization of the paper.} The paper is organized as follows. Section \ref{sec: Preliminaries} recalls preliminary results on Markov processes and random dynamical systems, along with several auxiliary lemmas. Section \ref{sec: Ergodicity of compact lifted processes} is devoted to establishing ergodicity of the associated base, one-point and projective processes. In Section \ref{sec: Quenched correlation decay: proof of Theorem 1.1}, we deduce exponential mixing from the exponential ergodicity of the
two-point process. Finally, Appendix \ref{sec: Proof of Lemma 2.13} contains the proof of the spanning property. 

\subsection*{Notation.} Throughout the paper, we denote by $\mathbb{S}^1$ the unit sphere in $\mathbb{R}^2$, and let $\|\cdot\|$ represent the Euclidean norm for vectors and the induced operator norm for matrices. With the regularization parameter $\epsilon>0$ fixed, we define
\begin{equation*}
    \Lambda_1:=2\|DK_\epsilon\|_{L^\infty},\qquad \Lambda_2:=2\left\|D^2K_\epsilon\right\|_{L^\infty}.
\end{equation*}Consequently, the first and second derivatives of $u(Z,x)$ are uniformly bounded by $\Lambda_1$ and $\Lambda_2$, respectively.

\section{Preliminaries}\label{sec: Preliminaries}
This section is devoted to several fundamental concepts from ergodic theory and random dynamical systems, as well as essential technical preliminaries. We first introduce the Markov processes associated with the randomized vortex system, which will play a central role in subsequent analysis.

\subsection{Associated Markov processes}\label{subsec: Associated Markov processes} Let $X_1(t),X_2(t)$ be the vortices defined in \eqref{eq: vortices}, and $\phi_t(x)$ the induced Lagrangian flow.

\medskip
\noindent\textbf{Base process.} The base process describes the joint evolution of the two randomized vortices, denoted by
\begin{equation*}
    Z_t=(X_1(t),X_2(t))\in\mathbb{T}^4.
\end{equation*}

\medskip
\noindent\textbf{One-point process.} The one-point process is defined as
\begin{equation*}
    (Z_t,x_t)\in\mathbb{T}^4\times\mathbb{T}^2, \qquad x_t:=\phi_t(x),
\end{equation*}where $x_t$ represents the position of the Lagrangian particle starting from $x\in\mathbb{T}^2$. We adjoin the vortex configuration $Z_t$ to the particle position $x_t$ to ensure that $(Z_t,x_t)$ constitutes a Markov process, whereas the trajectory $x_t$ alone is generally non-Markovian. For any initial state $(z,x)\in\mathbb{T}^4\times\mathbb{T}^2$ and Borel set $A\in\mathcal{B}\!\left(\mathbb{T}^4\times\mathbb{T}^2\right)$, the associated transition kernel is given by
\begin{equation*}\label{eq: def of one-point}
    P^{(1)}_t\!\left((Z,x),A\right)=\mathbb{P}\!\left((Z_t,x_t)\in A | (Z_0,x_0)=(Z,x)\right).
\end{equation*}

\medskip
\noindent\textbf{Projective process.} Given an initial unit tangent vector $v\in\mathbb{S}^1$ at $x\in\mathbb{T}^2$, define its evolution under the linearized flow by
\begin{equation*}
    v_t:=\frac{D_x\phi_t(x)v}{\|D_x\phi_t(x)v\|}.
\end{equation*}To study this directional dynamics, we introduce the projective process as the augmented process $(Z_t,x_t,v_t)$. Its transition kernel is explicitly given by
\begin{equation*}
    \hat{P}_t\!\left((Z,x,v),\hat{A}\right)=\mathbb{P}\!\left((Z_t,x_t,v_t)\in\hat{A}\bigl|(Z_0,x_0,v_0)=(Z,x,v)\right),
\end{equation*}where $\hat{A}\in \mathcal{B}\!\left(\mathbb{T}^4\times\mathbb{T}^2\times\mathbb{S}^1\right)$ is any Borel set.

\medskip
\noindent\textbf{Jacobian process.} We define the Jacobian process as 
\begin{equation*}
    (Z_t,x_t,J_t)\in\mathbb{T}^4\times\mathbb{T}^2\times \operatorname{SL}_2(\mathbb{R}),\qquad J_t:=D_x\phi_t(x).
\end{equation*}Due to the incompressibility of $u(t,x)$, the derivative matrix $J_t$ naturally evolves within $\operatorname{SL}_2(\mathbb{R})$. In contrast to the projective process, the Jacobian process retains both the directional and stretching information of the Lagrangian flow. This full structure is fundamental for rigorously studying the Lyapunov exponents. We denote the Markov transition kernel by $P_t^J$.

\medskip
\noindent\textbf{Two-point process.} This process traces the joint trajectories of two particles originating from distinct initial positions. More precisely, letting $\Delta:=\{(x,x):x\in\mathbb{T}^2\}$ be the diagonal in $\mathbb{T}^2\times\mathbb{T}^2$, the two-point process for initial positions $x\neq y$ is given by
\begin{equation*}
    (Z_t,x_t,y_t)\in\mathbb{T}^4\times\!\left((\mathbb{T}^2\times\mathbb{T}^2)\setminus\Delta\right),\qquad x_t:=\phi_t(x),\qquad y_t:=\phi_t(y).
\end{equation*}The corresponding transition kernel is denoted by $P_t^{(2)}$.

\vspace{0.3em}
Throughout this section, let $E$ denote a finite-dimensional complete metric space, $\mathcal{B}(E)$ its Borel $\sigma$-algebra, and $\mathcal{P}(E)$ the space of probability measures on $E$.

\subsection{Ergodicity and Harris theorem}\label{subsec: Ergodicity and Harris theorem} Assume that $\{P_t\}_{t\geq 0}$ is the family of transition kernels for a continuous-time Markov process on $E$. For any measurable observable $\varphi\in B_b(X)$ and probability measure $\mu\in\mathcal{P}(E)$, the associated Markov semigroup $P_t\colon B_b(E)\rightarrow B_b(E)$ and its dual $P_t^*\colon\mathcal{P}(E)\rightarrow\mathcal{P}(E)$ are respectively defined by
\begin{equation*}
    P_t\varphi(x)=\int_E\varphi(y)P_t(x,\mathrm{d}y),\qquad P_t^*\mu(A)=\int_EP_t(x,A)\mu(\mathrm{d}x).
\end{equation*}Throughout this section, we assume that $\{P_t\}_{t\geq 0}$ possesses the \textit{Feller property}, meaning that $P_tf\in C_b(E)$ for all $f\in C_b(E)$ and $t\geq 0$.

Recall that a probability measure $\mu\in\mathcal{P}(E)$ is said to be \textit{stationary} for $\{P_t\}_{t\geq 0}$ if $P_t^*\mu=\mu$ for any $t\geq 0$. A set $A\in\mathcal{B}(E)$ is $(P_t,\mu)$-\textit{invariant} if $P_t\mathbf{1}_A=\mathbf{1}_A$ holds $\mu$-almost surely for all $t\geq 0$. A stationary measure $\mu$ is then called \textit{ergodic} if every $(P_t,\mu)$-invariant set satisfies $\mu(A)=0$ or $\mu(A)=1$. Furthermore, given $V\colon E\rightarrow [0,+\infty)$, we define $L^{\infty}_V(E)$ as the space of all measurable functions $\varphi$ on $E$ for which the weighted supremum norm
\begin{equation}
    \|\varphi\|_V:=\sup_{x\in E}\frac{|\varphi(x)|}{1+V(x)}\label{eq: weighted norm}
\end{equation}is finite. The following theorem is a standard Harris-type result; for further details, we refer the reader to \cite{MT-93,HM-Harris}.
\begin{theorem}[Discrete Harris theorem]\label{thm: Discrete Harris theorem}
    Assume that there exists $T>0$ such that the transition kernel $P_T$ satisfies the following conditions:
    \begin{enumerate}[label=(\arabic*), ref=(\arabic*)]
        \item \label{assumption: Lyapunov drift}(Lyapunov drift condition) There exist a Lyapunov function $V\colon E\rightarrow [0,+\infty)$ and constants $K\geq 0,\gamma\in (0,1)$ such that
        \begin{equation}
            P_TV(x)\leq \gamma V(x)+K,\qquad \forall x\in E.
        \end{equation}
        \item \label{assumption: small set}(Minorization condition) There exists a constant $\alpha\in(0,1)$ and a probability measure $\nu$ so that
        \begin{equation}
            \inf_{x\in\mathcal{C}}P_T(x,\cdot)\geq \alpha\nu(\cdot),
        \end{equation}where $\mathcal{C}:=\{x\in E\colon V(x)\leq R\}$ for some $R>\frac{2K}{1-\gamma}$, and $K,\gamma$ are the constants from Assumption \ref{assumption: Lyapunov drift}.
    \end{enumerate}Then $P_T$ admits a unique invariant measure $\mu$. Furthermore, there exist constants $C>0$ and $\bar{\alpha}\in (0,1)$ such that for any $n\in\mathbb{N}$ and $\varphi\in L_V^\infty$, it holds that
    \begin{equation}
        \left\|P_{nT}\varphi-\int_{E}\varphi(x)\mu(\mathrm{d}x)\right\|_V\leq C\bar{\alpha}^n\left\|\varphi-\int_E\varphi(x)\mu(\mathrm{d}x)\right\|_V.
    \end{equation}  
\end{theorem}
\begin{rmk}
    We note that once the short-time estimate
    \begin{equation}\label{eq: abstract short-time estimate}
        \|P_t\varphi\|_V\leq C_T\|\varphi\|_V,\qquad \forall t\in [0,T)
    \end{equation}is established, the continuous-time contraction for all $t\in[0,+\infty)$ follows directly by combining Theorem \ref{thm: Discrete Harris theorem} with \eqref{eq: abstract short-time estimate}.
\end{rmk}
\begin{rmk}\label{rmk: SF+irreducibility}
    Typically, verifying the minorization condition is reduced to establishing the following properties:
    \begin{enumerate}[label=(\roman*), ref=(\roman*)]
        \item (\textbf{Strong Feller}) The Markov kernel $P_T$ is \textit{strong Feller} if for any $\psi\in B_b(E)$, the function $P_T\psi$ is continuous on $E$.
        \item (\textbf{Topological irreducibility}) For every $x\in X$ and any nonempty open set $A\subset E$, there exists $T>0$ such that $P_T(x,A)>0$.
    \end{enumerate} 
    We also mention that alternative ergodicity criteria are available under weaker regularity assumptions, including asymptotic strong Feller properties, eventual continuity, and suitable lower-bound conditions; see \cite{HM-06,HM-11EJP,GLLL-26,GLLL-25,LW-25}.
\end{rmk}
\begin{rmk}
    When the state space $E$ is compact, the Markov kernel $P_T$ automatically satisfies the Lyapunov-drift condition with $V\equiv1$. Consequently, the verification of Assumption \ref{assumption: Lyapunov drift} can be bypassed when establishing the ergodicity of the one-point and projective processes.
\end{rmk}
\subsection{Random dynamical systems}\label{subsec: Random dynamical systems} In this section, we introduce essential concepts from the theory of \textit{random dynamical systems} (RDS), which are crucial for establishing the ergodicity of the two-point process.

Let $(\Omega,\mathcal{F},\mathbb{P})$ be a probability space, where $\Omega$ is a Polish space and $\mathcal{F}=\mathcal{B}(\Omega)$. Furthermore, let 
\begin{equation*}
    \theta\colon\mathbb{R}_+\times\Omega\rightarrow\Omega,\qquad (t,\omega)\mapsto\theta_t\omega,
\end{equation*}be a measurable mapping such that
\begin{enumerate}[label=(\arabic*), ref=(\arabic*)]
    \item $\theta_0=\operatorname{Id}_\Omega$;
    \item $\mathbb{P}\circ(\theta_t)^{-1}=\mathbb{P}$ for all $t\geq 0$.
    \item $\theta_{t+s}=\theta_t\circ\theta_s$ for all $s,t\geq 0$.
\end{enumerate}The family $\{\theta_t\}_{t\geq 0}$ thus defines a measure-preserving semiflow on $\Omega$. With this setup in place, we now introduce the formal definition of a continuous RDS.
\begin{definition}\label{def: RDS}
    A continuous RDS on $E$ over $(\Omega,\mathcal{F},\mathbb{P},\{\theta_t\}_{t\geq 0})$ is a mapping
    \begin{equation*}
        \mathcal{A}\colon\mathbb{R}_+\times\Omega\times E\rightarrow E, \qquad (t,\omega,x)\mapsto\mathcal{A}_{\omega,t}x,
    \end{equation*}satisfying the following properties
    \begin{enumerate}[label=(\arabic*), ref=(\arabic*)]
        \item (Measurability) The mapping $\mathcal{A}$ is measurable with respect to $\mathcal{B}\!\left(\mathbb{R}_+\right)\otimes\mathcal{F}\otimes\mathcal{B}(E)$ and $\mathcal{B}(E)$.
        \item (Identity) $\mathcal{A}_{\omega,0}=\operatorname{Id}_E$ for all $\omega\in\Omega$.
        \item (Cocycle property) For all $s,t\geq 0$ and $\omega\in\Omega$,
        \begin{equation*}
            \mathcal{A}_{\omega,t+s}=\mathcal{A}_{\theta_s\omega,t}\circ\mathcal{A}_{\omega,s}
        \end{equation*}
        \item (Continuity) The mapping $\mathcal{A}_\omega\colon\mathbb{R}_+\times E\rightarrow E$ is continuous for all $\omega\in\Omega$.
    \end{enumerate}
\end{definition}
\begin{rmk}
    In addition to the conditions outlined in Definition \ref{def: RDS}, we assume throughout that the continuous RDS $\mathcal{A}$ possesses the independent increments property. That is, for all $s,t>0$, the random mapping $\mathcal{A}_{\cdot,t}$ is independent of $\mathcal{A}_{\theta_t\cdot,s}$.
\end{rmk}
For a continuous RDS $\mathcal{A}$ with independent increments and an initial state $x\in E$, we define the Markov process $x_t$ as follows (see, e.g., \cite{KS-12,BBPS-Lyapunov}): 
\begin{equation*}
    x_t(\omega):=\mathcal{A}_{\omega,t}x,\qquad x_0=x.
\end{equation*}Analogous to \eqref{eq: def of one-point}, for any $t>0,x\in E$ and $A\in \mathcal{B}(E)$, we consider the Markov kernel associated with $\mathcal{A}$
\begin{equation*}
    P_t(x,A)=\mathbb{P}(x_t\in A\mid x_0=x).
\end{equation*}The continuity of $\mathcal{A}$ immediately implies that $\{P_t\}_{t\geq 0}$ enjoys the Feller property. Furthermore, we say that $\mu\in \mathcal{P}(E)$ is \textit{stationary} (respectively, \textit{ergodic}) for the continuous RDS $\mathcal{A}$ if it is a stationary (respectively, ergodic) measure for the semigroup $\{P_t\}_{t\geq 0}$.

\vspace{0.3em}
Assuming that the mapping $x\mapsto\mathcal{A}_{\omega,t}x$ is of class $C^1$ on $E$, its derivative
\begin{equation*}
    D_x\mathcal{A}_{\omega,t}\colon \mathbb{R}_+\times\Omega\times E\rightarrow \operatorname{GL}_d(\mathbb{R}),\qquad d=\operatorname{dim}(E),
\end{equation*}is commonly referred to as the \textit{derivative cocycle}. Specifically, it satisfies the following \textit{cocycle property}:
\begin{equation*}
    D_x\mathcal{A}_{\omega,t+s}=D_{\mathcal{A}_{\omega,t}x}\mathcal{A}_{\theta_t\omega,s}\circ D_x\mathcal{A}_{\omega,t}.
\end{equation*}The long-time behavior of $D_x\mathcal{A}_{\omega,t}$ is described by the following theorem. For a comprehensive treatment, we refer the reader to \cite{Arnold-95}.
\begin{theorem}
    Let $\mathcal{A}$ be a continuous RDS and let $\mu\in\mathcal{P}(E)$ be an ergodic stationary measure associated with $\mathcal{A}$. Assume that the following integrability condition is satisfied:
    \begin{equation*}
        \int_{\Omega\times E}\left(\sup_{0\leq t\leq 1}\log^+\|D_x\mathcal{A}_{\omega,t}\|+\sup_{0\leq t\leq 1}\log^+\!\left\|(D_x\mathcal{A}_{\omega,t})^{-1}\right\|\right)\mathbb{P}(\mathrm{d}\omega)\mu(\mathrm{d}x)<+\infty.
    \end{equation*}Then there exist an integer $r$, $r$ distinct real numbers $\lambda_1>\cdots>\lambda_r>-\infty$, and for $\mathbb{P}\times\mu$-almost every $(\omega,x)\in\Omega\times E$, a measurable flag of subspaces
    \begin{equation*}
        \mathbb{R}^d=:F_1(\omega,x)\supset F_2(\omega,x)\supset\cdots\supset F_r(\omega,x)\supset F_{r+1}(\omega,x)=:\{0\}
    \end{equation*}such that the dimensions $\dim(F_i)$ are constant for $(\omega,x)$. Furthermore, for any $1\leq i\leq r$ and $v\in F_i(\omega,x)\setminus F_{i+1}(\omega,x)$, it holds that
    \begin{equation}
        \lambda_i=\lim_{t\rightarrow+\infty}\frac{1}{t}\log\!\left\|D_x\mathcal{A}_{\omega,t}v\right\|.\label{eq: definition of Lyapunov exponents} 
    \end{equation}In addition, the assignment $(\omega,x)\mapsto F_i(\omega,x)$ varies measurably.
\end{theorem}
The constants $\lambda_i$ are known as the \textit{Lyapunov exponents} with multiplicities
\begin{equation*}
    m_i:=\dim(F_i)-\dim(F_{i+1}).
\end{equation*}
Geometrically, these exponents characterize the asymptotic exponential growth rates of the singular values of $D_x\mathcal{A}_{\omega,t}$ as $t\rightarrow+\infty$. Consequently, we obtain the formula for the \textit{sum Lyapunov exponent}:
\begin{equation*}
    \lambda_\Sigma:=m_1\lambda_1+\cdots+m_r\lambda_r=\lim_{t\rightarrow+\infty}\frac{1}{t}\log\!\left|\det D_x\mathcal{A}_{\omega,t}\right|.
\end{equation*}If the continuous RDS $\mathcal{A}$ is volume-preserving, meaning that
\begin{equation*}
    \left|\det D_x\mathcal{A}_{\omega,t}\right|\equiv1,
\end{equation*}then the sum Lyapunov exponent vanishes identically. This immediately yields the bound
\begin{equation*}
    d\lambda_1\geq \lambda_\Sigma=0,
\end{equation*}where equality is achieved if and only if all Lyapunov exponents coincide. Furthermore, A strictly positive Lyapunov exponent ($\lambda_1>0$) reflects a sensitive dependence on initial conditions for the orbits, a property that is closely related to the mixing regime of the system; see, e.g., \cite{BBPS-Lyapunov,BBPS-mixing,BCZG-23}. To rule out the degenerate scenario where $\lambda_1=0$, Furstenberg's criterion was originally established for products of i.i.d. matrices \cite{Furstenberg-63}. Since then, this fundamental result has been extended to a variety of settings (see \cite{Ledrappier-86,Bougerol-88}). For the sake of simplicity, we adopt the formulation presented in \cite{BBPS-Lyapunov}, which is stated below. 

Considering the matrix process $(\mathcal{A}_{\omega,t}x,D_x\mathcal{A}_{\omega,t})\in E\times\operatorname{SL}_d(\mathbb{R})$ generated by the incompressible RDS $\mathcal{A}$, the corresponding transition kernel is given by
\begin{equation}
    Q_t\!\left((x,\operatorname{Id}),\cdot\right)=\mathbb{P}\!\left((\mathcal{A}_{\omega,t}x,D_x\mathcal{A}_{\omega,t})\in \cdot\right).
\end{equation}Specifically, fixing $t_0>0$, we define the regular conditional probability distribution as follows: 
\begin{equation}
    Q_{x,x'}^{t_0}(A):=\mathbb{P}\!\left(D_x\mathcal{A}_{\omega,t}\in A\mid \mathcal{A}_{\omega,t}x=x'\right),\qquad A\in\mathcal{B}(\operatorname{SL}_d(\mathbb{R})).
\end{equation}
\begin{theorem}[{\hspace{-0.05em}\cite[Corollary 3.22]{BBPS-Lyapunov}}]\label{thm: BBPS-condition} Let $\mathcal{A}$ be an incompressible continuous RDS with an ergodic stationary measure $\mu$. Assume that the derivative cocycle $D_x\mathcal{A}_{\omega,t}\in \operatorname{SL_d(\mathbb{R})}$ satisfies the following condition:
    \begin{enumerate}[label=(\roman*), ref=(\roman*)]
        \item\label{condition C} There exist a $t_0>0$ and a set $S\subset E$ with $\mu(S)>0$ such that for each $x\in E$, one can find a set $S_x\subset E$ with $P_{t_0}(x,S_x)>0$ having the property that for any $x'\in S_x$, the regular conditional probability measure $Q_{x,x'}^{t_0}(\cdot)$ is well-defined and absolutely continuous with respect to the Lebesgue measure on $\operatorname{SL}_d(\mathbb{R})$.
    \end{enumerate}Then we have $\lambda_1>0$.
\end{theorem}
\begin{rmk}
    To formalize the randomized vortex system within the dynamical framework, we identify the continuous RDS $\mathcal{A}$ with the one-point dynamics $(Z_t,x_t)$. The corresponding derivative cocycle with respect to the particle position is given by
    \begin{equation*}
        J_t=D_x\phi_t.
    \end{equation*}Accordingly, the associated matrix process has transition kernel $ Q_t=P_t^J$.
\end{rmk}

\subsection{Auxiliary tools}\label{subsec: Auxiliary tools} Recalling the Markov processes introduced in Section \ref{subsec: Associated Markov processes}, we can write the governing equations for $(x_t,v_t)\in\mathbb{T}^2\times\mathbb{S}^{1}$ as
\begin{align}
    \mathrm{d}x_t&=u(t,x_t)\mathrm{d}t,\label{eq: one-point dynamics}\\
    \mathrm{d}v_t&=\!\left(I-v_t\otimes v_t\right)D_xu(t,x_t)v_t\mathrm{d}t.\label{eq: projective dynamics}
\end{align}Likewise, the matrix process $J_t\in\operatorname{SL}_d(\mathbb{R})$ is given by
\begin{equation}
    \mathrm{d}J_t=D_xu(t,x_t)J_t\mathrm{d}t.\label{eq: Jacobian dynamics}
\end{equation}Therefore, the process of interest can be formulated as an abstract degenerate SDE
\begin{equation}\label{eq: general degenrate SDE}
    \mathrm{d}y_t=Y_0(y_t)\mathrm{d}t+\sum_{j=1}^4 Y_j\mathrm{d}B_j(t),
\end{equation}on the product space $\mathbb{T}^4\times \mathcal{M}$, where $\mathcal{M}$ is a finite-dimensional Riemann manifold (e.g., $\mathbb{T}^2$, $\mathbb{T}^2\times\mathbb{S}^1$ or $\mathbb{T}^2\times\operatorname{SL}_2(\mathbb{R})$). Here, $Y_0$ is the drift vector field on $\mathbb{T}^4\times\mathcal{M}$, while $\{Y_j\}_{j=1}^4$ constitutes a basis for the tangent space of $\mathbb{T}^4$, with their $\mathcal{M}$-components being identically zero.

Next, we recall the parabolic H\"ormander condition for the vector fields $\{Y_j\}$, which serves as a standard criterion for establishing the absolute continuity and the strong Feller property. For any two vector fields $Y$ and $Z$, we denote the standard Lie bracket by $[Y,Z]:=YZ-ZY$.
\begin{definition}
    Given a finite collection of vector fields $\{Y_j\}_{j=0}^n$ on a manifold $\mathcal{M}$, we define an increasing sequence of families of vector fields $\mathcal{Y}_0\subset\mathcal{Y}_1\cdots\subset\mathcal{Y}_k\subset\cdots$ recursively by
    \begin{align*}
        \mathcal{Y}_0&=\{Y_j:j\geq 1\},\\
        \mathcal{Y}_{k+1}&=\mathcal{Y}_k\cup\left\{[Y_j,Y]:Y\in\mathcal{Y}_{k},\ 0\leq j\leq n\right\}.
    \end{align*}We say that $\{Y_j\}_{j=0}^n$ satisfies the parabolic H\"ormander condition if for all $w\in \mathcal{M}$, one has
    \begin{equation*}
        \left\{Y(w):Y\in\bigcup_{k\geq0}\mathcal{Y}_k\right\}=T_w\mathcal{M}.
    \end{equation*}
\end{definition}
\begin{theorem}%[\hspace{-0.05em}\cite{Hormander-67}] 
Let $P_t(y,A)=\mathbb{P}(y_t\in A\mid y_0=y)$ be the Markov kernel associated with the finite-dimensional SDE \eqref{eq: general degenrate SDE}. If $\{Y_j\}_{j=0}^4$ is of class $C^\infty$ and satisfies the parabolic H\"ormander condition, then $P_t(y,\cdot)$ admits a smooth density with respect to $\operatorname{Leb}_{\mathbb{T}^4\times\mathcal{M}}$ for all $t>0$.
\end{theorem}
\begin{rmk}
    In the case where $\mathcal{M}$ is a compact manifold, the strong Feller property of $P_t$ follows immediately from the existence of a smooth transition density.
\end{rmk}
The following lemma plays a central role in verifying the H\"ormander condition for the randomized vortex system. The complete proof is deferred to Appendix.
\begin{lemma}\label{lem: nondegeneracy lemma}
    For any $r,s\in\mathbb{T}^2$ with $r\neq s$ and any $v\in\mathbb{S}^1$,the following identities hold:
    \begin{align}
        \operatorname{span}\!\left\{\partial^\alpha K_\epsilon(r):|\alpha|\geq 1\right\}&=\mathbb{R}^2,\label{eq:span one-point}\\
        \operatorname{span}\!\left\{\left(\partial^\alpha K_\epsilon(r),\Pi_v\partial^\alpha DK_{\epsilon}(r)v\right):|\alpha|\geq 1\right\}&=\mathbb{R}^2\oplus T_v\mathbb{S}^1,\label{eq:span projective}\\
        \operatorname{span}\!\left\{\left(\partial^\alpha K_\epsilon(r),\partial^\alpha DK_\epsilon(r)\right):|\alpha|\geq 1\right\}&=\mathbb{R}^2\oplus \mathfrak{sl}_2(\mathbb{R}),\label{eq:span Jocabian}\\
        \operatorname{span}\!\left\{\left(\partial^\alpha K_\epsilon(r),\partial^\alpha K_\epsilon(s)\right):|\alpha|\geq 1\right\}&=\mathbb{R}^2\oplus \mathbb{R}^2.\label{eq:span two-point}
    \end{align}Here $\Pi_v:=I-v\otimes v$ denotes the orthogonal projection onto $T_v\mathbb{S}^1$, and $\mathfrak{sl}_2(\mathbb{R}):=\{A\in M_2(\mathbb{R}):\operatorname{tr}(A)=0\}$ represents the special linear Lie algebra.
\end{lemma}

\section{Ergodicity of compact lifted processes}\label{sec: Ergodicity of compact lifted processes}
The primary goal of this section is to establish ergodicity for the base, one-point, and projective processes. Our strategy relies fundamentally on Harris's theorem. Since these processes evolve on compact state spaces, it suffices to verify the strong Feller property and topological irreducibility. We begin our analysis with the base process.
\subsection{The base process}\label{subsec: The base process} Let $m_4:=\operatorname{Leb}_{\mathbb{T}^4}$ denote the Lebesgue measure on $\mathbb{T}^4$. The following proposition identifies the stationary measure of the base process.
\begin{proposition}
    The base process $Z_t=(X_1(t),X_2(t))$ admits a unique stationary measure, which is precisely $m_4$.
\end{proposition}
\begin{proof}
    Set $Z=(X_1,X_2)\in\mathbb{T}^4$. The infinitesimal generator of the base process \eqref{eq: vortices} is given by
    \begin{equation*}
        \mathcal{L}_b=\frac{1}{2}\Delta_Z+b(Z)\cdot\nabla_Z.
    \end{equation*}Since the drift vector field $b$ is divergence-free, integration by parts yields
    \begin{equation*}
        \int_{\mathbb{T}^4}\mathcal{L}_bf\,\mathrm{d}m_4=0,
    \end{equation*}for any $f\in C^\infty$, which verifies that $m_4$ is stationary.

    To prove uniqueness, assume that $\mu$ is an arbitrary stationary measure ($P_t^*\mu=\mu$). Since $\mathcal{L}_b$ is uniformly elliptic with smooth coefficients, standard theory ensures that the transition kernel admits a positive smooth density with respect to $m_4$:
    \begin{equation*}
        P_t\!\left(Z,\mathrm{d}Z'\right)=p_t\!\left(Z,Z'\right)m_4\!\left(\mathrm{d}Z'\right).
    \end{equation*} The stationarity of $\mu$ therefore implies that
    \begin{equation*}
        \mu\!\left(\mathrm{d}Z'\right)=\left(\int_{\mathbb{T}^4}p_t\!\left(Z,Z'\right)\mu\!\left(\mathrm{d}Z\right)\right)m_4\!\left(\mathrm{d}Z'\right)=:\rho\!\left(Z'\right)m_4\!\left(\mathrm{d}Z'\right).
    \end{equation*}Furthermore, the stationary density $\rho$ satisfies the adjoint equation
    \begin{equation*}
         \mathcal{L}_b^*\rho
        =\frac12\Delta\rho-b\cdot\nabla\rho=0.
    \end{equation*}Multiplying this equation by $\rho$ and integrating over 
    $\mathbb{T}^4$ gives
    \begin{equation*}
        -\frac{1}{2}\int_{\mathbb{T}^4}\!\left|\nabla\rho\right|^2\mathrm{d}m_4=0.
    \end{equation*}Consequently, $\rho\equiv c$ for some constant $c$. Given that $\rho$ is a probability density on $\mathbb{T}^4$, we conclude $\rho\equiv 1$, and hence $\mu=m_4$.
\end{proof}
\begin{rmk}
    In fact, the conclusion can be significantly strengthened. By combining uniform ellipticity, compactness of the state space, and positivity of the heat kernel, it follows that the base process is uniformly exponentially ergodic in total variation distance:
    \begin{equation*}
        \sup_{Z\in\mathbb{T}^4}\|P_t(Z,\cdot)-m_4(\cdot)\|_{\operatorname{TV}}\lesssim e^{-\gamma t}.
    \end{equation*}
\end{rmk}

\subsection{The one-point process}\label{subsec: The one-point process}
We now turn to the study of the one-point Lagrangian motion $(Z_t,x_t)$. The proof of ergodicity consists of two main steps. First, we verify the parabolic H\"ormander condition, which implies the strong Feller property. Second, utilizing a controllability argument, we establish topological irreducibility of the one-point process. Combining these two properties, we conclude the ergodicity of the one-point process.

\subsubsection{Parabolic H\"ormander condition} Recall that for an initial state $x\in\mathbb{T}^2$, the trajectory $x_t$ is defined by
\begin{equation*}
    \mathrm{d}x_t=u(Z_t,x)\mathrm{d}t.
\end{equation*}Therefore, the drift vector field associated with $(Z_t,x_t)$ is
\begin{equation*}
    F_1=b(Z)\cdot\nabla_Z+u(Z,x)\cdot\nabla_x.
\end{equation*}Furthermore, we denote by $G_{i,a}:=\partial_{X_{i,a}}$ (for $i,a\in \{1,2\}$) the four constant noise vector fields.
\begin{proposition}\label{prop: one-point parabolic Hormander}
    The collection of vector fields $\{F_1, G_{i,a}:i,a=1,2\}$ satisfies the parabolic H\"ormander condition.
\end{proposition}
\begin{proof}
    Since the vector fields $G_{i,a}$ span all directions corresponding to the vortex variables $(X_1,X_2)$, it remains to verify that the $x$-components of the resulting iterated commutators span $T_x\mathbb{T}^2$. A direct computation gives
    \begin{equation*}
        \Pi_x\left[G_{i,a},F_1\right]=-\partial_{a}K_\epsilon(x-X_i),
    \end{equation*}where $\Pi_x$ denotes the projection onto the $x$-component. Repeated Lie bracketing with $G_{i,1}$ and $G_{i,2}$ vector fields of the form 
    \begin{equation*}
        \partial^{\alpha}K_\epsilon(x-X_i),\qquad |\alpha|\geq 1.
    \end{equation*}By \eqref{eq:span one-point}, this family of fields spans $T_x\mathbb{T}^2$. Consequently, the parabolic H\"ormander condition holds.
\end{proof}

\subsubsection{Topological irreducibility} In this section, we establish the exact controllability of the one-point process, which subsequently implies the irreducibility of the transition kernel $P_t^{(1)}$. More precisely, we demonstrate that by constructing suitable control paths, the system can be steered from any initial configuration to any prescribed target state. Our construction relies on the observation that
\begin{equation*}
    \operatorname{span}\left\{K_\epsilon(r):r\in\mathbb{T}^2\right\}=\mathbb{R}^2.
\end{equation*}Therefore,  there exist two points $\xi_1,\xi_2\in\mathbb{T}^2$ such that the vectors
\begin{equation}
    \eta_1=2K_\epsilon(\xi_1),\qquad \eta_2=2K_\epsilon(\xi_2),
\end{equation}constitute a basis for $\mathbb{R}^2$.

\begin{proposition}\label{prop: one-point controllability}
    The controlled one-point system $Y=(Z,x)$:
    \begin{equation}\label{eq: controlled one-point}
        \begin{cases}
            \dot{X_1}=K_\epsilon\!\left(X_1-X_2\right)+h_1,\\
            \dot{X_2}=K_\epsilon\!\left(X_2-X_1\right)+h_2,\\
            \ \ \dot{x}=u(Z,x), 
        \end{cases}
    \end{equation}is exactly controllable. Specifically, for every pair of states 
    \begin{equation*}
        Y^0=\left(X_1^0,X_2^0,x^0\right),\qquad Y^*=\left(X_1^*,X_2^*,x^*\right)\in\mathbb{T}^4\times\mathbb{T}^2,
    \end{equation*}there exist $T>0$ and controls $h_1,h_2\in L^2\!\left([0,T];\mathbb{R}^2\right)$ such that the corresponding solution to \eqref{eq: controlled one-point} with $Y_0=Y^0$ satisfies $Y_T=Y^*$.
\end{proposition}
\begin{proof}
    Let
    \begin{equation*}
        r_i^0:=x^0-X_i^0,\qquad r_i^*=x^*-X_i^*,
    \end{equation*}denote the initial and target relative coordinates. We first construct a path for the relative positions, and then recover the vortex trajectories and the corresponding controls from it. 

    Consider a continuous, piecewise $C^1$ path 
    \begin{equation*}
        \bar{r}\colon [0,T_0]\rightarrow\mathbb{T}^4,
    \end{equation*}that connects $\left(r_1^0,r_2^0\right)$ and $\left(r_1^*,r_2^*\right)$ by passing successively through
    \begin{equation}
        \left(r_1^0,r_2^0\right)\rightarrow (\xi_1,\xi_1)\rightarrow -(\xi_1,\xi_1)\rightarrow (\xi_2,\xi_2)\rightarrow -(\xi_2,\xi_2)\rightarrow \!\left(r_1^*,r_2^*\right).
    \end{equation}Integrating along this path yields a displacement in the $x$-component:
    \begin{equation*}
        C:=\int_0^{T_0}\left(K_\epsilon(\bar r_1(t))+K_\epsilon(\bar r_2(t))\right)\mathrm{d}t.
    \end{equation*}Fix lifts $\tilde x^0,\tilde x^*\in\mathbb{R}^2$ of $x^0,x^*$ and let
    \begin{equation*}
        D:=\tilde{x}^*-\tilde{x}^0.
    \end{equation*}Since $\eta_1,\eta_2$ form a basis of $\mathbb{R}^2$, there exist $\alpha,\beta\in\mathbb{R}$ such that
    \begin{equation*}
        D-C=\alpha\eta_1+\beta\eta_2.
    \end{equation*}Furthermore, we decompose
    \begin{equation*}
        \alpha=\alpha_+-\alpha_-,\qquad \beta=\beta_+-\beta_-,
    \end{equation*}with $\alpha_{\pm},\beta_{\pm}\geq 0$.

    To compensate for this gap $D-C$, we modify $\bar{r}$ by inserting the following constant waiting intervals:
    \begin{equation}
         \begin{array}{c|c}
            \text{relative state} & \text{waiting time}\\ \hline
            (\xi_1,\xi_1)       & \alpha_+\\
            -(\xi_1,\xi_1)      & \alpha_-\\
            (\xi_2,\xi_2)       & \beta_+\\
            -(\xi_2,\xi_2)      & \beta_-
        \end{array}.
    \end{equation}Let
    \begin{equation*}
        r=(r_1,r_2)\colon[0,T]\rightarrow \mathbb{T}^4,
    \end{equation*}denote the resulting path, where $T=T_0+\alpha_++\alpha_-+\beta_++\beta_-$. The endpoints remain
    \begin{equation*}
        r(0)=\!\left(r_1^0,r_2^0\right),\qquad r(T)=\!\left(r_1^*,r_2^*\right).
    \end{equation*}Recalling the odd symmetry of $K_\epsilon$ (which implies $2K_\epsilon(-\xi_i)=-\eta_i$), we obtain
    \begin{equation*}
        \int_0^T \left(K_\epsilon(r_1(t))+K_\epsilon(r_2(t))\right)\mathrm{d}t=C+\alpha_+\eta_1-\alpha_-\eta_1+\beta_+\eta_2-\beta_-\eta_2=D.
    \end{equation*}Define a lifted particle trajectory by
    \begin{equation*}
        \tilde{x}_t=\tilde{x}^0+\int_0^t \left(K_\epsilon(r_1(s))+K_\epsilon(r_2(s))\right)\mathrm{d}s.
    \end{equation*}Its restriction to $\mathbb{T}^2$, denoted by $x_t$, trivially satisfies
    \begin{equation*}
        \dot{x}=K_\epsilon(r_1)+K_\epsilon(r_2),\quad x_0=x^0,\quad x_T=x^*.
    \end{equation*}
    
    Finally, by setting $X_i(t):=x_t-r_i(t)$, we can construct the control paths
    \begin{align}
        h_1&:=K_\epsilon(r_1)+K_\epsilon(r_2)-K_\epsilon(r_2-r_1)-\dot{r}_1,\label{eq: one-point control h1}\\
        h_2&:=K_\epsilon(r_1)+K_\epsilon(r_2)-K_\epsilon(r_1-r_2)-\dot{r}_2.\label{eq: one-point control h2}
    \end{align}It then follows that
    \begin{align*}
        \dot{X_1}&=\dot{x}-\dot{r_1}=K_\epsilon(r_2-r_1)+h_1=K_\epsilon(X_1-X_2)+h_1,\\
        \dot{X_2}&=\dot{x}-\dot{r_2}=K_\epsilon(r_1-r_2)+h_2=K_\epsilon(X_2-X_1)+h_2.
    \end{align*}Evaluating at the endpoints yields $X_i(0)=X_i^0$ and $X_i(T)=X_i^*$. Moreover, the piecewise $C^1$ regularity of $r(t)$ guarantees that $h_i\in L^2$. Consequently, the controlled trajectory $(X_1,X_2,x)$ steers $\left(X_1^0,X_2^0,x^0\right)$ to $\left(X_1^*,X_2^*,x^*\right)$.
\end{proof}
Observe that the drift field of the one-point process is divergence-free. Following the same argument as in Section \ref{subsec: The base process}, we deduce that the Lebesgue measure $m_6:=\operatorname{Leb}_{\mathbb{T}^4\times\mathbb{T}^2}$ is stationary. Combining this observation with Propositions \ref{prop: one-point parabolic Hormander} and \ref{prop: one-point controllability}, we obtain the following immediate corollary.
\begin{corollary}
    The one-point process $(Z_t,x_t)$ admits $m_6$ as its unique stationary measure.
\end{corollary}

\subsection{The projective process}\label{subsec: The projective process}
This section is devoted to establishing the ergodicity of the projective process. Our argument parallels the one employed in Section \ref{subsec: The one-point process}: it suffices to verify the parabolic H\"ormander condition and the irreducibility. To the end, we begin with the parabolic H\"ormander condition.
\subsubsection{Parabolic H\"ormander condition}
    The projective dynamics on $\mathbb{T}^2\times\mathbb{S}^1$ are governed by \eqref{eq: one-point dynamics} and \eqref{eq: projective dynamics}. Accordingly, for $Z=(X_1,X_2)$, the associated drift vector field is explicitly given by
    \begin{equation}
        F_P=b(Z)\cdot\nabla_Z+u(Z,x)\cdot\nabla_x+\Pi_v\!\left(D_xu(Z,x)v\right)\cdot\nabla_v.
    \end{equation}
    \begin{proposition}\label{prop: projective parabolic Hormander}
        The drift vector field $F_P$, together with the noise fields $G_{i,a}$, satisfies the parabolic H\"ormander condition.
    \end{proposition}
    \begin{proof}
        As in the proof of Proposition \ref{prop: one-point parabolic Hormander}, we can restrict our attention to the $(x,v)$-components. By taking iterated Lie brackets of $F_P$ with $G_{i,a}$, we obtain the following $(x,v)$-projections:
        \begin{equation*}
            \left(\partial^{\alpha}K_\epsilon(x-X_i),\Pi_v\partial^\alpha DK_\epsilon(x-X_i)v\right), \qquad |\alpha|\geq 1.
        \end{equation*}Invoking Lemma \ref{lem: nondegeneracy lemma}, we conclude that this family of vectors spans the joint tangent space $T_x\mathbb{T}^2\oplus T_v\mathbb{S}^1$.        
    \end{proof}
\subsubsection{Topological irreducibility} To prove the irreducibility of the projective process, we consider the associated deterministic control system $Y=(Z,x,v)$:
\begin{equation}\label{eq: controlled projective}
    \begin{cases}
            \dot{X_1}=K_\epsilon\!\left(X_1-X_2\right)+h_1,\\
            \dot{X_2}=K_\epsilon\!\left(X_2-X_1\right)+h_2,\\
            \ \ \dot{x}=u(Z,x),\\
            \ \ \dot{v}=\Pi_vD_xu(Z,x)v.
    \end{cases}
\end{equation}The following proposition demonstrates the exact controllability of this system.
\begin{proposition}\label{prop: projective controllability}
    For any $Y^0,Y^*\in\mathbb{T}^4\times\mathbb{T}^2\times\mathbb{S}^1$, there exist a time $T>0$ and controls $h_1,h_2\in L^2\!\left([0,T];\mathbb{R}^2\right)$ such that the trajectory of the controlled system \eqref{eq: controlled projective} starting from $Y^0$ satisfies $Y_T=Y^*$. In particular, the Markov kernel $\hat{P}_T$ is topologically irreducible.
\end{proposition}
The following result is a key ingredient in the argument.
\begin{lemma}\label{lem: coincident rotation}
    Let 
    \begin{equation*}
        J=\begin{pmatrix}
            0 & 1\\
            -1 & 0
        \end{pmatrix}.
    \end{equation*}Then there exists a constant $\gamma_\epsilon>0$ such that
    \begin{equation*}
        DK_\epsilon(0)=-\gamma_\epsilon J.
    \end{equation*}Therefore, at the configuration $X_1=X_2=x$, the dynamics simplify to $\dot{x}=0$ and 
    \begin{equation}
        \dot{v}=-2\gamma_\epsilon Jv,\qquad  v_t=e^{-2\gamma_\epsilon tJ}v_0.\label{eq: rotation of v_t}
    \end{equation}As $t\geq 0$ varies, the family $\left\{e^{-2\gamma_\epsilon tJ}\right\}_{t\geq 0}$ coincides with the full rotation group on $\mathbb{S}^1$.
\end{lemma}
\begin{proof}[Proof of Lemma \ref{lem: coincident rotation}]
    Differentiating the Fourier series at the origin gives
    \begin{equation*}
        DK_\epsilon(0)=\begin{pmatrix}
            \sum_{k\neq\mathbf{0}}e^{-4\pi^2\epsilon^2|k|^2}\frac{k_1k_2}{|k|^2} & -\sum_{k\neq\mathbf{0}}e^{-4\pi^2\epsilon^2|k|^2}\frac{k_1^2}{|k|^2}\\
            \sum_{k\neq\mathbf{0}}e^{-4\pi^2\epsilon^2|k|^2}\frac{k_2^2}{|k|^2} & -\sum_{k\neq\mathbf{0}}e^{-4\pi^2\epsilon^2|k|^2}\frac{k_1k_2}{|k|^2}
        \end{pmatrix}.
    \end{equation*}By symmetry under the reflection $k_1\rightarrow -k_1$ and the interchange $k_1\rightarrow k_2$, we deduce that the two diagonal sums vanish and the off-diagonal entries are opposite. Setting
    \begin{equation*}
        \gamma_\epsilon:=\sum_{k\neq \mathbf{0}}e^{-4\pi^2\epsilon^2|k|^2}\frac{k_1^2}{|k|^2}>0,
    \end{equation*}the Jacobian matrix reduces to
    \begin{equation*}
        DK_\epsilon(0)=-\gamma_\epsilon J.
    \end{equation*}Since $Jv\perp v$ for any $v\in\mathbb{S}^1$, the projection $\Pi_v$ preserves \(Jv\), which yields \eqref{eq: rotation of v_t}.
\end{proof}
\begin{proof}[Proof of Proposition \ref{prop: projective controllability}]
    Let 
    \begin{equation*}
        Y^0=\!\left(X_1^0,X_2^0,x^0,v^0\right),\qquad Y^*=\!\left(X_1^*,X_2^*,x^*,v^*\right),
    \end{equation*}be the initial and target states. Fix any $\bar{x}\in\mathbb{T}^2$. By Proposition \ref{prop: one-point controllability}, there exist $T_1>0$ and controls $h_1^{(1)},h_2^{(1)}\in L^2\!\left([0,T_1];\mathbb{R}^2\right)$ such that the one-point trajectory of \eqref{eq: controlled one-point} reaches $(\bar{x},\bar{x},\bar{x})$ at time $T_1$. We denote the corresponding projective direction by $v_{T_1}$.
    
    To resolve the target state, we apply a similar argument to find a time $T_2>0$ and controls $h_1^{(2)},h_2^{(2)}$ that steer the one-point state from $(\bar{x},\bar{x},\bar{x})$ to $\left(X_1^*,X_2^*,x^*\right)$. Fix one such pair of controls. The corresponding linearized flow induces a diffeomorphism
    \begin{equation*}
        \Phi_2\colon\mathbb{S}^1\rightarrow\mathbb{S}^1,\qquad v\mapsto\frac{J_{T_2}v}{\|J_{T_2}v\|}.
    \end{equation*}Therefore, there exists a unique
    \begin{equation*}
        w:=\Phi_2^{-1}(v^*),
    \end{equation*}such that if the second stage initiates with direction $w$, then the terminal projective direction is exactly $v^*$.

    It remains to bridge the gap between $v_{T_1}$ and $w$ while preserving the spatial position. To achieve this, upon reaching $(\bar{x},\bar{x},\bar{x})$, we set $h_1=h_2=0$. At this coincident configuration, 
    \begin{equation*}
        \dot{X_1}=\dot{X_2}=\dot{x}=0,
    \end{equation*}while the direction evolves according to
    \begin{equation*}
        \dot{v}=-2\gamma_\epsilon J v,\qquad v_t=e^{-2\gamma_\epsilon tJ}v_{T_1}.
    \end{equation*}By Lemma \ref{lem: coincident rotation}, we can select a time $T_0>0$ such that
    \begin{equation*}
        e^{-2\gamma_\epsilon T_0 J}v_{T_1}=w.
    \end{equation*}

    We now concatenate the three control stages. Let $T:=T_1+T_0+T_2$ and define, for $i=1,2,$
    \begin{equation*}
        h_i(t):= 
        \begin{cases}
            h_i^{(1)}(t), &t\in [0,T_1),\\
            0, &t\in[T_1,T_1+T_0],\\
            h_i^{(2)}(t-T_1-T_0) &t\in(T_1+T_0,T].
        \end{cases}
    \end{equation*}These controls effectively steer the full state $Y^0$ exactly to $Y^*$, completing the proof.   
\end{proof}
\begin{corollary}
    The projective process has a unique invariant probability measure, denoted by $\pi_P$.
\end{corollary}
\vspace{0.3em}
To connect the derivative growth with the projective dynamics, we introduce the instantaneous stretching observable
\begin{equation*}
    \sigma\!\left(Z,x,v\right)=v\cdot D_xu(Z,x)v.
\end{equation*}
Since the Jacobian $J_t$ is governed by
\begin{equation*}
    J_t=D_x\phi_t,
\end{equation*}a direct computation yields
\begin{equation*}
    \frac{\mathrm{d}}{\mathrm{d}t}\log\|J_tv\|=\sigma(Z_t,x_t,v_t),\qquad \log\|J_tv\|=\int_0^t\sigma(Z_s,x_s,v_s)\mathrm{d}s.
\end{equation*}The top Lyapunov exponent $\lambda_1$ can therefore be expressed as the stationary average of $\sigma$, as stated in the following Furstenberg--Khasminskii formula.
\begin{theorem}[{\hspace{-0.05em}\cite[Proposition 2.1]{BBPS-ICM}}]\label{thm: Furstenberg--Khasminskii formula.}
    For $m_6$-a.e. $(Z,x)$ and all $v\in\mathbb{S}^1$, there holds 
    \begin{equation*}
        \lambda_1=\lim_{t\rightarrow\infty}\frac{1}{t}\log\|J_tv\|\qquad \text{almost surely.}
    \end{equation*}Furthermore, we have
    \begin{equation}
        \lambda_1=\int_{\mathbb{T}^4\times\mathbb{T}^2\times\mathbb{S}^1}\sigma\,\mathrm{d}\pi_P.\label{eq: Furstenberg--Khasminskii formula}
    \end{equation}
\end{theorem}

\section{Quenched correlation decay: proof of Theorem \ref{thm: exponential mixing}}\label{sec: Quenched correlation decay: proof of Theorem 1.1}
We now turn to the proof of Theorem \ref{thm: exponential mixing}. A crucial intermediate step is to establish ergodicity of the two-point process. For any pair of distinct initial points $(x,y)\in\Delta^c=(\mathbb{T}^2\times\mathbb{T}^2)\setminus\Delta$, let
\begin{equation*}
    x_t=\phi_t(x),\qquad y_t=\phi_t(y),
\end{equation*}denote the two-point trajectory. The associated lifted process $\left(Z_t,x_t,y_t\right)$ induces the Feller Markov semigroup
\begin{equation*}           
    P_t^{(2)}\varphi\left(Z,x,y\right):=\mathbb{E}_{Z,x,y}\varphi\!\left(Z_t,x_t,y_t\right),\quad \varphi\in B_b\!\left(\mathbb{T}^4\times\Delta^c\right).
\end{equation*}
The ergodicity argument consists of two main steps. First, in Section \ref{subsec: Drift condition for the two-point process}, we verify the drift condition (Assumption \ref{assumption: Lyapunov drift}) by constructing an explicit Lyapunov function. Second, Section \ref{subsec: Fixed-time minorization} is devoted to establishing the minorization condition (Assumption \ref{assumption: small set}) at a fixed time $T$. Applying Theorem \ref{thm: Discrete Harris theorem} to these estimates yields geometric ergodicity of the time-$T$ skeleton $P_{nT}^{(2)}$. We then extend this to all continuous time $t\geq 0$, which ultimately yields the desired exponential mixing.
\subsection{Drift condition for the two-point process}\label{subsec: Drift condition for the two-point process} A primary challenge is controlling the two-point process near the diagonal $\Delta$. To this end, we first establish positivity of the top Lyapunov exponent and then utilize this positivity to derive a uniform contraction estimate.
\subsubsection{Positivity of Lyapunov exponents}\label{sec: Positivity of Lyapunov exponents}
We prove that the top Lyapunov exponent is strictly positive by applying Theorem \ref{thm: BBPS-condition}. Since the one-point
process admits a unique stationary measure $m_6$ on $\mathbb{T}^4\times\mathbb{T}^2$ (see Section \ref{subsec: The one-point process}), the problem reduces to verifying Condition \ref{condition C}. Our analysis begins by establishing the hypoellipticity of the Jacobian process.

Based on the evolution equation \eqref{eq: Jacobian dynamics}, the drift vector field for the Jacobian process is explicitly formulated as
\begin{equation*}
    F_J=b(Z)\cdot\nabla_Z+u(Z,x)\cdot\nabla_x+D_xu(Z,x)J\cdot\nabla_J.
\end{equation*}For any $J\in\operatorname{SL}_2(\mathbb{R})$, we have
\begin{equation*}
    T_J\operatorname{SL}_2(\mathbb{R})=\mathfrak{sl}_2(\mathbb{R})J.
\end{equation*}Under the right trivialization $HJ\mapsto H$, the tangent space is identified with the Lie algebra $\mathfrak{sl}_2(\mathbb{R})$. Exploiting this algebraic structure yields the following result.
\begin{proposition}\label{prop: Jacobian Hormander}
    The Jacobian process satisfies the parabolic H\"ormander condition.
\end{proposition}
\begin{proof}
    Modulo the vortex noise directions, the iterated Lie brackets of $F_J$ with $G_{i,a}$ generate $(x,J)$-components of the form
    \begin{equation*}
        \left(\partial^{\alpha}K_\epsilon(x-X_i),\partial^{\alpha}DK_{\epsilon}(x-X_i)J\right),\qquad |\alpha|\geq 1.
    \end{equation*}After applying the right trivialization, this family of vector fields is represented as 
    \begin{equation*}
        \left(\partial^{\alpha}K_\epsilon(x-X_i),\partial^{\alpha}DK_{\epsilon}(x-X_i)\right).
    \end{equation*}Since these fields span $T_x\mathbb{T}^2\oplus \mathfrak{sl}_2(\mathbb{R})$, we conclude that the parabolic H\"ormander condition holds.
\end{proof}

With the H\"ormander condition secured for the Jacobian process, we now proceed to verify Condition \ref{condition C}. The crucial observation is that hypoellipticity yields a smooth transition density for $(Z_t,x_t,J_t)$. By disintegrating this joint law with respect to the one-point marginal, we obtain the conditional absolute continuity. The precise statement is as follows.
\begin{proposition}
    The derivative cocycle $J_t=D_x\phi_t$ satisfies Condition \ref{condition C} with $t_0=1$. In particular, we can set $S=\mathbb{T}^4\times\mathbb{T}^2$ such that for every $(Z,x)\in S$, the subset $S_{(Z,x)}$ satisfies
    \begin{equation*}
        P^{(1)}_1\!\left((Z,x),S_{(Z,x)}\right)=1.    
    \end{equation*}
\end{proposition}
\begin{proof}
    By Proposition \ref{prop: Jacobian Hormander} and the finite-dimensional H\"ormander theorem (see, e.g., \cite{Hairer-Hormander}), the time-one transition kernel of the Jacobian process admits the representation
    \begin{equation*}                    
        P_1^J\!\left((Z,x,I),\mathrm{d}Z'\mathrm{d}x'\mathrm{d}A'\right)=q_{Z,x}\!\left(Z',x',A'\right)m_6\!\left(\mathrm{d}Z'\mathrm{d}x'\right)\mathrm{d}A'.
    \end{equation*}Integrating out the matrix variable yields the one-point marginal measure
    \begin{equation*}
        P_1^{(1)}\!\left((Z,x),\mathrm{d}Z'\mathrm{d}x'\right)=p_{Z,x}\!\left(Z',x'\right)m_6\!\left(\mathrm{d}Z'\mathrm{d}x'\right),
    \end{equation*}where $p_{Z,x}(\cdot)=\int_{\operatorname{SL}_2(\mathbb{R})}q_{Z,x}\!\left(\cdot,A'\right)\mathrm{d}A'$. Since $P_1^{(1)}$ is a probability measure, $p_{Z,x}$ is finite $m_6$-almost everywhere. Define 
    \begin{equation*}
        S_{Z,x}:=\left\{(Z',x')\in\mathbb{T}^4\times\mathbb{T}^2:0<p_{Z,x}\!\left(Z',x'\right)<\infty\right\}.
    \end{equation*}It naturally follows that $S_{Z,x}$ is a set of full measure, i.e., $P_1^{(1)}\!\left((Z,x),S_{Z,x}\right)=1$. Furthermore, disintegrating the joint measure yields the conditional probability
    \begin{equation*}
         Q_{(Z,x),(Z',x')}^1\!\left(\mathrm{d}A'\right)=\frac{q_{Z,x}\!\left(Z',x',A'\right)}{p_{Z,x}\!\left(Z',x'\right)}\mathrm{d}A'.
    \end{equation*}Therefore,
    \begin{equation*}
        Q^1_{(Z,x),(Z',x')}\ll\operatorname{Leb}_{\operatorname{SL}_2(\mathbb{R})},
    \end{equation*}for every $\left(Z',x'\right)\in S_{Z,x}$, which completes the proof.
\end{proof}

\subsubsection{Lyapunov drift condition}
We now turn to the construction of the Lyapunov function. The following uniform contraction condition is a key ingredient in the argument.
\begin{lemma}\label{lem: uniform contraction condition}
    There exist a time $T>0$ and constants $p,\gamma'\in \!\left(0,\frac{1}{2}\right)$ such that uniform inverse-moment estimate holds
    \begin{equation}
        \mathbb{E}_{Z,x}\|D_x\phi_T(x)v\|^{-p}\leq \gamma',\label{eq: uniform contraction}
    \end{equation}for all $(Z,x)\in\mathbb{T}^4\times\mathbb{T}^2$ and $v\in\mathbb{S}^1$.
\end{lemma}
\begin{proof}
    Combining the Furstenberg--Khasminskii formula \eqref{eq: Furstenberg--Khasminskii formula} with standard Ces\`aro convergence, we obtain
    \begin{equation*}
        \lim_{t\rightarrow \infty}\sup_{(Z,x,v)}\left|\frac{1}{t}\mathbb{E}_{Z,x}\log\|J_tv\|-\lambda_1\right|=\lim_{t\rightarrow\infty}\sup_{(Z,x,v)}\left|\frac{1}{t}\int_0^t\sigma(Z_s,x_s,v_s)\mathrm{d}s-\int\sigma\,\mathrm{d}\pi_P\right|=0.
    \end{equation*}Thus, there exists $T_0>0$ such that all initial projective states $(Z,x,v)$,
    \begin{equation}\label{eq: lower bound for J_Tv}
        \mathbb{E}_{Z,x}\log\|J_{T_0}v\|\geq \frac{\lambda_1T_0}{2}.
    \end{equation}

    The boundedness of the stretching observable $\sigma$ implies
    \begin{equation*}
        \big|\log\|J_{T_0}v\|\big|\leq\Lambda_1T_0.
    \end{equation*}By fixing $p>0$ such that $p\Lambda_1T_0<1$, Taylor's theorem gives
    \begin{equation}\label{eq: upper bound for J_Tv}
        \|J_{T_0}v\|^{-p}=e^{-p\log\|J_tv\|}\leq 1-p\log\|J_{T_0}v\|+\frac{e}{2}p^2\Lambda_1^2T_0^2.
    \end{equation}We further restrict $p$ to be small enough so that the second-order error term is absorbed
    \begin{equation*}
        \frac{e}{2}p^2\Lambda_1^2T_0^2\leq \frac{p\lambda_1T_0}{4}<1.
    \end{equation*}Taking the expectation of \eqref{eq: upper bound for J_Tv} and applying \eqref{eq: lower bound for J_Tv}, we deduce the uniform one-step contraction
    \begin{equation*}
        \mathbb{E}_{Z,x}\|J_{T_0}v\|^{-p}\leq 1-\frac{p\lambda_1T_0}{4}=:\gamma_0<1, 
    \end{equation*}To complete the proof, choosing $n\in\mathbb{N}$ such that $\gamma_0^n<1/2$ and setting $T=nT_0$, an application of the Markov property yields
     \begin{equation*}
        \mathbb{E}\!\left[\|J_{(k+1)T_0}v\|^{-p}\mid \mathcal{F}_{kT_0}\right]\leq \gamma_0\|J_{kT_0}v\|^{-p}. 
    \end{equation*}Iterating this estimate concludes the proof of \eqref{eq: uniform contraction} with $\gamma':=\gamma_0^n$.
\end{proof}
We next establish the drift condition for the two-point process. To control the singular behavior near the diagonal, the uniform contraction condition \eqref{eq: uniform contraction} suggests considering 
\begin{equation*}
    V(Z,x,y):=\operatorname{dist}_{\mathbb{T}^2}(x,y)^{-p},
\end{equation*}where the exponent $p$ is chosen as in Lemma \ref{lem: uniform contraction condition}. This particular choice yields the following drift-type estimate.
\begin{proposition}\label{prop: global drift}
    There exist $T>0$, $\gamma\in (0,1)$ and $C_0>0$ such that 
    \begin{equation}
        P_T^{(2)}V\leq \gamma V+C.\label{eq: global drift condition}
    \end{equation}     
\end{proposition}
\begin{proof}
    We adopt the constants $T,p,\gamma'$ given by Lemma \ref{lem: uniform contraction condition}. For any $x\neq y$, choose lifts $\tilde{x},\tilde{y}\in\mathbb{R}^2$ satisfying 
    \begin{equation*}
        \|\tilde{y}-\tilde{x}\|=\operatorname{dist}_{\mathbb{T}^2}(x,y),
    \end{equation*}and let $\tilde{x_t},\tilde{y}_t$ denote the corresponding lifted trajectories. We track their separation $\delta_t:=\tilde{y}_t-\tilde{x}_t$. Applying Taylor's expansion yields
    \begin{equation*}
        \dot{\delta_t}=D_xu(Z_t,x_t)\delta_t+R_t,\qquad \|R_t\|\leq \frac{\Lambda_2}{2}\|\delta_t\|^2.\label{eq: delta_t equation}
    \end{equation*}Integrating \eqref{eq: delta_t equation} and utilizing Lipschitz estimates $\|\delta_t\|\leq e^{\Lambda_1t}\|\delta_0\|$, we obtain
    \begin{align}
        \left\|\delta_T-J_T\delta_0\right\|&\leq \int_0^T \left\|J_TJ_s^{-1}\right\|\cdot\|R_s\|\mathrm{d}s\notag\\[0.2em]
        &\leq \frac{\Lambda_2T}{2}e^{2\Lambda_1T}\|\delta_0\|^2\leq C_*\|\delta_0\|^2,\label{eq: upper bound for delta_T-J_T}       
    \end{align}where we set $C_*:=2+\frac{\Lambda_2T}{2}e^{2\Lambda_1T}$. 

    We then choose $r_*>0$ sufficiently small so that
    \begin{equation*}
        r_*\leq \frac{e^{-\Lambda_1T}}{2C_*}.
    \end{equation*}For any initial points with $\operatorname{dist}_{\mathbb{T}^2}(x,y)=\|\delta_0\|\leq r_*$, the separation remains small over the interval $[0,T]:$
    \begin{equation*}
        \|\delta_t\|\leq e^{\Lambda_1t}r_*\leq \frac{1}{4}.
    \end{equation*}This ensures that the lifted and torus distances agree on $[0,T]$. Note that $\left\|J_T\delta_0\right\|\geq e^{-\Lambda_1T}\|\delta_0\|$, we deduce that
    \begin{align}
        \left\|\delta_T\right\|\geq \left\|J_T\delta_0\right\|-C_*\|\delta_0\|^2&\geq \left\|J_T\delta_0\right\|-C_*r_*\|\delta_0\|\notag \\[0.2em]
        &\geq \left\|J_T\delta_0\right\|-\frac{1}{2}e^{-\Lambda_1T}\|\delta_0\|\notag\\
        &\geq \frac{1}{2}\left\|J_T\delta_0\right\|,\label{eq: lower bound for delta_t}
    \end{align}By normalizing the initial displacement to $v=\frac{\delta_0}{\|\delta_0\|}\in\mathbb{S}^1$, and combining \eqref{eq: lower bound for delta_t} with Lemma \ref{lem: uniform contraction condition}, we conclude that
    \begin{align}
        P_T^{(2)}V(Z,x,y)=\mathbb{E}\!\left\|\delta_T\right\|^{-p}&\leq 2^{p}\mathbb{E}\left\|J_Tv\right\|^{-p}\|\delta_0\|^{-p}\notag\\[0.2em]
        &\leq 2^p\gamma'\operatorname{dist}_{\mathbb{T}^2}(x,y)^{-p}\notag\\[0.2em]
        &=\gamma V(Z,x,y),\label{eq: upper bound near the diagonal}
    \end{align}where $\gamma:=2^p\gamma'\in (0,1)$.

    Away from the diagonal, for $x,y$ with 
    \begin{equation*}
        \operatorname{dist}{\mathbb{T}^2}(x,y)\geq r_*,
    \end{equation*}the inverse Lipschitz bound gives pathwise
    \begin{equation}
        V(Z_T,x_T,y_T)\leq e^{p\Lambda_1T}r_*^{-p}=:C_0.\label{eq: upper bound away from diagonal}
    \end{equation}Combining \eqref{eq: upper bound near the diagonal} and \eqref{eq: upper bound away from diagonal}, we obtain the global drift condition
    \begin{equation*}
        P_T^{(2)}V\leq \gamma V+C_0.
    \end{equation*}
\end{proof}
\subsection{Fixed-time minorization}\label{subsec: Fixed-time minorization}
The purpose of this section is to verify the minorization condition in Theorem \ref{thm: Discrete Harris theorem}. Following Remark \ref{rmk: SF+irreducibility}, the argument is reduced to establishing the strong Feller property and topological irreducibility of the two-point process. These are proved in Sections \ref{subsubsec: Strong Feller property} and \ref{subsubsec: Topological irreducibility}, respectively. We then combine these two ingredients to derive a fixed-time minorization estimate.
\subsubsection{Strong Feller property}\label{subsubsec: Strong Feller property} To prove the strong Feller property, we begin by verifying the parabolic H\"ormander condition. The drift vector field associated with the two-point process is given by
\begin{equation*}
    F_2:=b(Z)\cdot\nabla_Z+u(Z,x)\cdot\nabla_x+u(Z,y)\cdot\nabla_y. 
\end{equation*}
\begin{lemma}
    The drift field $F_2$, together with the diffusion vector fields $G_{i,a}$, satisfies the parabolic H\"ormander condition.
\end{lemma}
\begin{proof}
    As in the preceding H\"ormander arguments, the noise fields $G_{i,a}$ already span all vortex directions. It therefore suffices to focus on the passive-particle directions. By computing repeated Lie brackets, we obtain $(x,y)$-components
    \begin{equation*}
        \left(D^{\alpha}K_\epsilon(x-X_i),D^{\alpha}K_{\epsilon}(y-X_i)\right),\qquad |\alpha|\geq 1.
    \end{equation*}As $x\neq y$, \eqref{eq:span two-point} guarantees that these vectors span $\mathbb{R}^2\oplus\mathbb{R}^2$, which completes the proof.
\end{proof}
Since the system features smooth coefficients, the H\"ormander condition yields a smooth transition density, which is far stronger than mere absolute continuity. Although the state space $\Delta^c$ is noncompact, the Lipschitz estimates for the Lagrangian flow provide the additional control needed to establish the strong Feller property.
\begin{proposition}\label{prop: two-point strong Feller}
    For all $t>0$, the Markov kernel $P_t^{(2)}$ admits a smooth density
    \begin{equation*}
        p_t\!\left(\zeta,\zeta'\right)\in C^\infty\!\left((0,\infty)\times\!\left(\mathbb{T}^4\times\Delta^c\right)^2\right)
    \end{equation*}and satisfies the strong Feller property.
\end{proposition}
\begin{proof}
    Fix a time $t>0$ and an initial state $\zeta_0=(Z_0,x_0,y_0)\in\mathbb{T}^4\times\Delta^c$, and set 
    \begin{equation*}
        d_0:=\operatorname{dist}_{\mathbb{T}^2}(x_0,y_0)>0.
    \end{equation*}By the parabolic H\"ormander theorem \cite{Hormander-67}; see also \cite{Stroock-08}, $P_t^{(2)}$ admits a smooth density
    \begin{equation*}
        P_t^{(2)}(\zeta_0,\mathrm{d}\zeta)=p_t\!\left(\zeta_0,\zeta\right)\mathrm{d}z,\qquad \forall \zeta\in\mathbb{T}^4\times\Delta^c.
    \end{equation*}

    Next, we choose a compact neighborhood $U\subset\Delta^c$ of $(x_0,y_0)$ such that
    \begin{equation*}
        \operatorname{dist}_{\mathbb{T}^2}(x,y)\geq \frac{d_0}{2},\qquad \forall (x,y)\in U.
    \end{equation*}The inverse Lipschitz estimate the yields
    \begin{equation*}
        \operatorname{dist}_{\mathbb{T}^2}(x_t,y_t)\geq e^{-\Lambda_1t}\operatorname{dist}_{\mathbb{T}^2}(x_0,y_0).
    \end{equation*}Hence for any trajectory starting from $\mathbb{T}^4\times U$, the state at time $t$ is strictly supported in the compact subset $\mathbb{T}^4\times K_\delta$, where
    \begin{equation*}
        K_\delta:=\!\left\{(x,y)\in\Delta^c:\operatorname{dist}_{\mathbb{T}^2}(x,y)\geq \delta\right\},\qquad \delta:=\frac{d_0}{2}e^{-\Lambda_1t}.
    \end{equation*}

    Let $f\in B_b\!\left(\mathbb{T}^4\times\Delta^c\right)$ and $\zeta_n\rightarrow \zeta_0$. For $n$ sufficiently large, $\zeta_n\in \mathbb{T}^4\times U$, which allows us to bound the difference by
      \begin{align*}
        \left|P_t^{(2)}f(\zeta_n)-P_t^{(2)}f(\zeta_0)\right|&\leq\int_{\mathbb{T}^4\times K_\delta} |f\!\left(\zeta\right)|\cdot\left|p_t\!\left(\zeta_n,\zeta\right)-p_t\!\left(\zeta_0,\zeta\right)\right|\mathrm{d}\zeta\\[0.2em]
        &\leq\|f\|_\infty\,\operatorname{Leb}\!\left(\mathbb{T}^4\times K_\delta\right)\sup_{\zeta\in \mathbb{T}^4\times K_\delta}\left|p_t\!\left(\zeta_n,\zeta\right)-p_t\!\left(\zeta_0,\zeta\right)\right|.
    \end{align*}By smoothness of $p_t$ and compactness, the supremum converges to zero as $n\rightarrow \infty$. This proves the continuity of $P_t^{(2)}f$ at $z_0$, thereby establishing the strong Feller property.
\end{proof}

\subsubsection{Topological irreducibility}\label{subsubsec: Topological irreducibility} Following standard methodology, we establish topological irreducibility of the two-point process via an approximate control argument. The associated controlled system is given by
\begin{equation}\label{eq: controlled two-point}
    \begin{cases}
            \dot{X_1}=K_\epsilon\!\left(X_1-X_2\right)+h_1,\\
            \dot{X_2}=K_\epsilon\!\left(X_2-X_1\right)+h_2,\\
            \ \ \dot{x}=u(Z,x),\\
            \ \ \dot{y}=u(Z,y).
    \end{cases}
\end{equation}

We begin the proof with the following auxiliary lemma. For any $z\in\mathbb{T}^2\setminus\{\mathbf{0}\}$, define
\begin{equation*}
    \mathcal{V}_z=\left\{
        \begin{pmatrix}
            K_\epsilon(a)\\
            K_\epsilon(a+z)
        \end{pmatrix}
    :a\in\mathbb T^2\right\}\subset\mathbb R^4.
\end{equation*}Let $\mathcal{W}_z=\mathcal{V}_z+\mathcal{V}_z$ denote its Minkowski sum; equivalently,
\begin{equation*}
    \mathcal{W}_z=\left\{W_z(a,b):a,b\in\mathbb T^2\right\}, \qquad W_z(a,b)=\begin{pmatrix}
        K_\epsilon(a)+K_\epsilon(b)\\
        K_\epsilon(a+z)+K_\epsilon(b+z)
    \end{pmatrix}.
\end{equation*}
\begin{lemma}\label{lem: Uniform convex velocity ball}
    For any $z\neq \mathbf{0}$, 
    \begin{equation*}
        \operatorname{span}\mathcal{V}_z=\mathbb{R}^4,\qquad \mathbf{0}\in\operatorname{int}\operatorname{co}(\mathcal{V}_z),\qquad \mathbf{0}\in\operatorname{int}\operatorname{co}(\mathcal{W}_z),
    \end{equation*}where $\operatorname{co}(\cdot)$ denotes the convex hull. Moreover, for every compact subset $E\subset\mathbb{T}^2\setminus\{\mathbf{0}\}$, there exists $\rho_E>0$ such that
    \begin{equation*}
        B\!\left(\mathbf{0},\rho_E\right)\subset\operatorname{co}(\mathcal{W}_z),\qquad \forall z\in E.
    \end{equation*}
\end{lemma}
\begin{proof}
    Fix $z\neq\mathbf{0}$. Suppose there exist $l,m\in\mathbb{R}^2$ satisfying
    \begin{equation*}
        l\cdot K_\epsilon(a)+m\cdot K_\epsilon(a+z)=0,\qquad \forall a\in\mathbb{T}^2.
    \end{equation*}
    Taking Fourier coefficients then yields
    \begin{equation*}
        l\cdot k^{\perp}+e^{2\pi ik\cdot z}m\cdot k^{\perp}=0,\qquad k\neq\mathbf{0}.
    \end{equation*}Following the argument in the proof of Lemma \ref{lem: nondegeneracy lemma}, we can select two linearly independent integer vectors $k_1,k_2$ such that
    \begin{equation*}
        l\cdot k_j^{\perp}=m\cdot k_j^{\perp}=0,\qquad j=1,2.
    \end{equation*}This implies that $l=m=\mathbf{0}$, thus confirming that $\mathcal{V}_z$ spans $\mathbb{R}^4$.
    
    By the odd symmetry of $K_\epsilon$, we have
    \begin{equation*}
        \int_{\mathbb{T}^2}\begin{pmatrix}
            K(a)\\
            K(a+z)
        \end{pmatrix}\mathrm{d}a=\mathbf{0},
    \end{equation*}which ensures that $\mathbf{0}\in\operatorname{co}(\mathcal{V}_z)$. If the origin were not an interior point, the supporting hyperplane theorem would yield a nonzero linear functional $H$ such that 
    \begin{equation*}
        H(V)\geq 0,\qquad \forall V\in\mathcal{V}_z.
    \end{equation*}Since the integral of $H(V)$ over $\mathbb{T}^2$ is zero, $H(V)$ must vanish identically on $\mathcal{V}_z$. This contradicts the fact that $\operatorname{span}\mathcal{V}_z=\mathbb{R}^4$, and hence $\mathbf{0}\in\operatorname{int}\operatorname{co}(\mathcal{V}_z)$.

    Given that $\mathcal{W}_z=\mathcal{V}_z+\mathcal{V}_z$, the origin is an interior point of $\operatorname{co}(\mathcal{W}_z)$ as well. Finally, the continuity of the mapping $z\mapsto\mathcal{W}_z$ in the Hausdorff metric, combined with the compactness of $E\subset\mathbb{T}^2\setminus\{\mathbf{0}\}$, guarantees the existence of a uniform radius $\rho_E>0$.
\end{proof}
\begin{rmk}
    In the two-point system, the variables $a,b$ and $z$ may be interpreted as the relative displacements
    \begin{equation*}
        y-x,\quad x-X_1,\quad x-X_2,
    \end{equation*}respectively.
\end{rmk}
The preceding lemma reveals that, whenever the two passive particles are separated by a positive distance, the vortex positions can be rapidly adjusted so as to realize arbitrary joint velocities. Building upon this observation, we now establish the approximate controllability of the two-point process.
\begin{proposition}\label{prop: two-point control}
    Given any initial state $\zeta^0\in\mathbb{T}^4\times\Delta^c$ and any nonempty open set $O\subset\mathbb{T}^4\times\Delta^c$, there exist a time $T>0$ and controls $h_1,h_2\in L^2\!\left([0,T];\mathbb{R}^2\right)$ such that the controlled trajectory \eqref{eq: controlled two-point} originating from $\zeta^0$ satisfies 
    \begin{equation*}
        (Z_T,x_T,y_T)\in O.
    \end{equation*}In particular, the two-point process is topologically irreducible.
\end{proposition}
\begin{proof}
    Write $\zeta^0=\!\left(X_1^0,X_2^0,x^0,y^0\right)$ and choose a target state $\zeta^*=\!\left(X_1^*,X_2^*,x^*,y^*\right)\in O$. We denote the initial and target relative configurations by
    \begin{equation*}
        r_i^0=x^0-X_i^0,\qquad r_i^*=x^*-X_i^*,\qquad i=1,2.
    \end{equation*}By openness of $O$, there exists $\eta>0$ such that any state $(x,y)$ satisfying 
    \begin{equation*}
        \operatorname{dist}_{\mathbb{T}^2}(x,x^*)+\operatorname{dist}_{\mathbb{T}^2}(y,y^*)<\eta,
    \end{equation*}implies
    \begin{equation*}
        \left(x-r_1^*,x-r_2^*,x,y\right)\in O.
    \end{equation*}

    We begin by selecting a smooth reference path 
    \begin{equation*}
        \bar{q}(t)=\!\left(\bar{x}(t),\bar{y}(t)\right)\in\Delta^c,\qquad t\in[0,T_0], 
    \end{equation*}connecting $\left(x^0,y^0\right)$ to $\left(x^*,y^*\right)$. Define the relative distance $\bar{z}(t):=\bar{y}(t)-\bar{x}(t)$. Since its image
    \begin{equation*}
         E=\{\bar{z}(t):t\in [0,T_0]\}
    \end{equation*}is compact in $\mathbb{T}^2\setminus\{\mathbf{0}\}$, Lemma \ref{lem: Uniform convex velocity ball} gives $\rho_E>0$ such that
    \begin{equation*}
        B(\mathbf{0},\rho_E)\subset \operatorname{co}\!\left(\mathcal{W}_{\bar{z}(t)}\right)
    \end{equation*}for every $t$. Reparametrizing $\bar q$ at a sufficiently slow speed, we may assume that 
    \begin{equation*}
        \left|\dot{\bar{q}}(t)\right|<\frac{\rho_E}{2},
    \end{equation*}which implies $\dot{\bar{q}}(t)\in \operatorname{co}\!\left(W_{\bar{z}(t)}\right)$. 
    
    Consider a partition $0=t_0<\cdots<t_N=T_0$ with mesh size $\Pi:=\max_n\{t_{n+1}-t_n\}$. According to Carath\'eodory's theorem, the velocity vector $\dot{\bar{q}}(t_n)$ can be expressed as a convex combination of at most five elements from $\mathcal{W}_{\bar{z}(t_n)}$:
    \begin{equation*}
        \dot{\bar{q}}(t_n)=\sum_{j=1}^{J_n}\lambda_{n,j}W_{\bar{z}(t_n)}(a_{n,j},b_{n,j}),
    \end{equation*}where $J_n\leq 5$, $\lambda_{n,j}\geq 0$ and $\sum_{j}\lambda_{n,j}=1$. On $[t_n,t_{n+1}]$, we define piecewise constant relative controls $\left(\hat{r}_1,\hat{r}_2\right)=(a_{n,j},b_{n,j})$ for a duration of $\lambda_{n,j}(t_{n+1}-t_n)$. This yields a trajectory $\hat{q}=(\hat{x},\hat{y})$ governed by
    \begin{equation*}
        \frac{\mathrm{d}}{\mathrm{d}t}\begin{pmatrix}
            \hat{x}\\
            \hat{y}
        \end{pmatrix}=W_{\hat{y}-\hat{x}}\!\left(\hat{r}_1(t),\hat{r}_2(t)\right).
    \end{equation*}Given that this vector field is Lipschitz continuous and uniformly bounded, standard Euler estimates and Gronwall's inequality ensure that
    \begin{equation*}
        \sup_{0\leq t\leq T_0}\operatorname{dist}\!\left(\hat{q}(t),\bar{q}(t)\right)\rightarrow 0,\qquad \operatorname{as}\ \Pi\rightarrow0.
    \end{equation*}

    It remains to mollify the jumps and prescribe the boundary values. Let
    \begin{equation*}
        c_j=(a_j,b_j),\qquad j=1,\ldots, M,
    \end{equation*}be chronological sequence of values taken by the piecewise control $\hat{q}$, and set
    \begin{equation*}
        c_0=\!\left(r_1^0,r_2^0\right),\qquad c_{N+1}=\!\left(r_1^*,r_2^*\right).
    \end{equation*}For each $0\leq j\leq M$, we connect $c_j$ to $c_{j+1}$ by a smooth path 
    \begin{equation*}
        \gamma_j\colon[0,\epsilon_j]\rightarrow \mathbb{T}^4,
    \end{equation*}where the duration $\epsilon_j$ can be chosen sufficiently small. By inserting these short transition paths $\gamma_j$ into $\hat{q}$, we obtain a piecewise smooth control $r(t)=(r_1(t),r_2(t))\in H^1\!\left([0,T];\mathbb{T}^4\right)$ satisfying 
    \begin{equation*}
        r(0)=\!\left(r_1^0,r_2^0\right),\qquad r(T)=\!\left(r_1^*,r_2^*\right),
    \end{equation*}where $T=T_0+\sum_j\epsilon_j$. 
    
    Finally, letting $(x_t,y_t)$ be the solution driven by $r(t)$
    \begin{equation*}
        \frac{\mathrm{d}}{\mathrm{d}t}\begin{pmatrix}
            x\\
            y
        \end{pmatrix}=W_{y-x}(r_1(t),r_2(t)),
    \end{equation*}and setting $X_i(t):=x_t-r_i(t)$ for $i=1,2$, the explicit controls are defined as
    \begin{align*}
        h_1&=K_\epsilon(r_1)+K_\epsilon(r_2)-K_\epsilon(r_2-r_1)-\dot{r}_1,\\
        h_2&=K_\epsilon(r_1)+K_\epsilon(r_2)-K_\epsilon(r_1-r_2)-\dot{r}_2.
    \end{align*}When the transition times $\epsilon_j$ are sufficiently small, Gronwall's lemma confirms that 
    \begin{equation}
        \operatorname{dist}_{\mathbb{T}^2}(x_T,x^*)+\operatorname{dist}_{\mathbb{T}^2}(y_T,y^*)<\eta.\label{eq: distance between q and q*}
    \end{equation}Combining \eqref{eq: distance between q and q*} with the choice of $\eta$, we deduce that
    \begin{equation*}
        \left(X_1(T),X_2(T),x_T,y_T\right)=\!\left(x_T-r_1^*,x_T-r_2^*,x_T,y_T\right)\in O.
    \end{equation*}Thus, the desired approximate controllability follows.
\end{proof}
\subsubsection{Minorization condition}\label{subsubsec: Minorization condition} We are now prepared to verify Assumption \ref{assumption: small set} of Theorem \ref{thm: Discrete Harris theorem}. The argument proceeds in two steps. We first establish a localized minorization condition. By subsequently combining this local estimate with the strong Feller property and topological irreducibility, we derive the global minorization condition for a prescribed skeleton step.
\begin{lemma}\label{lem: local minorization}
    There exist a compact set $B\subset\mathbb{T}^4\times\Delta^c$ with nonempty interior, a closed time interval $[a,b]\subset (0,\infty)$, a constant $\epsilon_B\in (0,1)$ and a probability measure $\nu_B$ such that
    \begin{equation}\label{eq: local minorization}
        P_t^{(2)}\!\left(\zeta,\cdot\right)\geq \epsilon_B\nu_B(\cdot),\qquad t\in [a,b],\quad \zeta\in B.
    \end{equation}
\end{lemma}
\begin{proof}
    Fix $\zeta_0\in\mathbb{T}^4\times\Delta^c$ and $s>0$. Since $p_s(\zeta_0,\cdot)$ is a smooth probability density, there exists a state $\zeta_1\in\mathbb{T}^4\times\Delta^c$ such that
    \begin{equation*}
        p_s(\zeta_0,\zeta_1)>0.
    \end{equation*}By the continuity of $p_s$, we can find an open neighborhood $U$ of $z_0$ satisfying
    \begin{equation*}
        p_s(\zeta,\zeta_1)>0,\qquad \forall \zeta\in U.
    \end{equation*}Then, topological irreducibility of the two-point process yields $r>0$ with
    \begin{equation*}
        P_r^{(2)}(\zeta_1,U)>0,
    \end{equation*}which implies $p_r(\zeta_1,\zeta_*)>0$ for some $\zeta_*\in U$. Given that $p_s(\zeta_*, \zeta_1) > 0$ as well, continuity and the Chapman--Kolmogorov identity guarantees that
    \begin{align}
        p_{r+s}(\zeta_*,\zeta_*)&=\int_{\mathbb{T}^4\times\Delta^c}p_s(\zeta_*,\zeta)p_r(\zeta,\zeta_*)\mathrm{d}\zeta\notag\\[0.2em]
        &\geq \int_{B(\zeta_1,r)}p_s(\zeta_*,\zeta)p_r(\zeta,\zeta_*)\mathrm{d}\zeta>0.\label{eq: positive-return}
    \end{align}
    
    Furthermore, joint smoothness from Proposition \ref{prop: two-point strong Feller} gives an open neighborhood $O$ of $z_*$, a time interval $[a,b]$ containing $ r+s$, and a constant $c>0$ such that 
    \begin{equation}
        p_t(\zeta,\zeta')\geq c,\qquad \zeta,\zeta'\in O,\quad t\in[a,b].\label{eq: positive local-density}
    \end{equation}Choosing a sufficiently small closed ball $B \subset O$, we define 
    \begin{equation*}
        \nu_B(\cdot):=\frac{\operatorname{Leb}(B\cap\cdot)}{\operatorname{Leb}(B)},\qquad \epsilon_B:=\min\left\{\frac{1}{2},c\operatorname{Leb}(B)\right\}.
    \end{equation*}Clearly, $B$ is compact with a nonempty interior and satisfies $\operatorname{supp} \nu_B \subset B$. Integrating \eqref{eq: positive local-density} over $A\cap B$, we obtain 
    \begin{equation*}
        P_t^{(2)}(\zeta,A)\geq P_t^{(2)}(\zeta,A\cap B)\geq \epsilon_B\nu_B(A),
    \end{equation*}thereby establishing \eqref{eq: local minorization}.
\end{proof}
The following result is a direct corollary of Lemma \ref{lem: local minorization}.
\begin{lemma}\label{lem: long-time local minorization}
    Let $B,\nu_B$ be given by Lemma \ref{lem: local minorization}. There exists $L_B>0$ such that for any $t\geq L_B$, we can find a constant $\epsilon_t\in (0,1)$ for which
    \begin{equation}
        P_t^{(2)}(\zeta,\cdot)\geq \epsilon_t\nu_B(\cdot),\qquad \forall\zeta\in B.\label{eq: long-time local minorization}
    \end{equation}
\end{lemma}
\begin{proof}
    Since $\nu_B$ is supported on $B$, repeatedly applying the Chapman--Kolmogorov identity to \eqref{eq: local minorization} yields
    \begin{equation*}
        P_t^{(2)}(\zeta,\cdot)\geq \epsilon_B^k\nu_B(\cdot),\qquad t\in [ka,kb],\quad \zeta\in B,
    \end{equation*}for every $k\in\mathbb{Z}_+$. 
    
    We can choose a sufficiently large integer $k_0 \geq 1$ such that $k_0b\geq (k_0+1)a$. This ensures that 
    \begin{equation*}
        kb\geq (k+1)a,\qquad \forall k\geq k_0,
    \end{equation*}which implies the intervals $[ka,kb]$ overlap, yielding 
    \begin{equation*}
        \bigcup_{k\ge k_0}[ka,kb]=[k_0a,\infty).
    \end{equation*}Setting $L_B:=k_0a$, it follows that for any $t\geq L_B$, there exists some $k\in\mathbb{Z}_+$ satisfying $t\in[ka,kb]$. Therefore, we obtain
    \begin{equation*}
        P_t^{(2)}(\zeta,\cdot)\geq \epsilon_B^k\nu_B(\cdot)=:\epsilon_t\nu_B(\cdot),
    \end{equation*}which completes the proof.
\end{proof}

Combining the preceding minorization estimates with strong Feller property and topological irreducibility, we now establish the global minorization condition at any prescribed skeleton
time step.
\begin{proposition}\label{prop: global minorization}
    Given any nonempty compact set $E\subset \mathbb{T}^4\times\Delta^c$ and any time $\tau>0$, there exist an integer $m>0$ and $\epsilon_E\in (0,1)$ such that 
    \begin{equation*}
        P_{m\tau}^{(2)}(\zeta,\cdot)\geq \epsilon_E\nu_B(\cdot),\qquad \zeta\in E.\label{eq: global minorization}
    \end{equation*}
\end{proposition}
\begin{proof}
    For each $\zeta\in E$, Proposition \ref{prop: two-point control} gives a time $t_\zeta>0$ for which
    \begin{equation*}
        P_{t_\zeta}^{(2)}(\zeta,B)>0.
    \end{equation*}Combining strong Feller property of $P_t^{(2)}$ with the compactness of $E$, we extract a finite open cover $U_1,\ldots,U_N$ of $E$, associated times $t_j > 0$ and constants $p_j$ such that
    \begin{equation*}
        P_{t_j}^{(2)}(\zeta,B)\geq p_j,\qquad \forall \zeta\in U_j,\quad 1\leq j\leq N.
    \end{equation*}
    
    Next, we select an integer $m$ large enough that
    \begin{equation*}
        m\tau\ge L_B+\max_j\{t_j\},
    \end{equation*}and define the residual times $s_j:=m\tau-t_j$. If $\zeta\in E\cap U_j$, applying the Chapman--Kolmogorov identity and the local bound from Lemma \ref{lem: long-time local minorization} then yields
    \begin{align*}
        P_{m\tau}^{(2)}(\zeta,A)&\geq \int_B P_{s_j}^{(2)}\!\left(\zeta',A\right)P_{t_j}^{(2)}\!\left(\zeta,\mathrm{d}\zeta'\right)\\
        &\geq p_j\epsilon_{s_j}\nu_B(A).
    \end{align*}Taking $\epsilon_E:= \min_j \{p_j \epsilon_{s_j}\}$ completes the proof.
\end{proof}

\subsection{Proof of Theorem \ref{thm: exponential mixing}}\label{subsec: Proof of Theorem 1.1}
We conclude this section by completing the proof of Theorem \ref{thm: exponential mixing}. Throughout the proof, $C$ denotes a generic constant, independent of $K$, whose value may change from line to line.
\begin{proof}
    The proof proceeds in two main steps. First, we apply Harris theorem to a suitable skeleton kernel $P_T^{(2)}$ and subsequently extend this result to all continuous times $t\geq 0$. Second, following the standard strategy initiated in \cite{DKK-04} and applied in recent literature \cite{BBPS-mixing,BCZG-23,CIS-25,LS-26}, we establish the desired quenched exponential mixing estimate.
    \begin{enumerate}[label=\textit{Step (\arabic*).}, ref=(\arabic*), wide=0pt, labelindent=\parindent, listparindent=\parindent]
        \item\label{step 1: proof of Theorem 1.1} Let $p,T,\gamma,C_0$ be the constants from Proposition \ref{prop: global drift}. Iterating \eqref{eq: global drift condition} yields, for any integer $m\geq 1$, 
        \begin{equation}
            P_{mT}^{(2)}V\leq \gamma^mV+\frac{C_0\!\left(1-\gamma^m\right)}{1-\gamma}.\label{eq: drift condition for mT}
        \end{equation}Choose $R>\max\left\{2^p,\frac{2C_0}{1-\gamma}\right\}$ and consider the compact level set
        \begin{equation*}
            E_R:=\!\left\{V\leq R\right\}=\!\left\{(Z,x,y):\operatorname{dist}_{\mathbb{T}^2}(x,y)\geq  R^{-1/p}\right\}.
        \end{equation*}Applying Proposition \ref{prop: global minorization} on $E=E_R$ with $\tau=T$ gives an integer $m\geq 1$, a constant $\epsilon_R\in (0,1)$ and a probability measure $\nu_R$ satisfying
        \begin{equation}
            P_{mT}^{(2)}(\zeta,\cdot)\geq \epsilon_R\nu_R(\cdot),\qquad \zeta\in E_R.\label{eq: global minorization condition for mT}
        \end{equation}Setting $T_*=mT$, \eqref{eq: drift condition for mT} and \eqref{eq: global minorization condition for mT} show that the fixed kernel $P_{T_*}^{(2)}$ satisfies the hypotheses of the discrete Harris theorem (Theorem \ref{thm: Discrete Harris theorem}). Since the product Lebesgue measure is invariant under the two-point semigroup, uniqueness identifies the ergodic measure $\pi$ with $\operatorname{Leb}$. Consequently, for any mean-zero observable $f\in H^s$, there exist constants $C$ and $\alpha\in (0,1)$ such that
        \begin{equation}
            \left\|P_{nT_*}^{(2)}f\right\|_V\leq C\alpha^n\|f\|_V,\label{eq: drift condition for nT*}
        \end{equation}where $\|\cdot\|_V$ denotes the weighted norm \eqref{eq: weighted norm}.

        To generalize \eqref{eq: drift condition for nT*} to continuous time $t\geq 0$, observe that for any $\zeta\in\mathbb{T}^4\times\Delta^c$,
        \begin{align}
            \left|\frac{P_t^{(2)}f\!\left(\zeta\right)}{1+V(\zeta)}\right|&\leq \int\left|\frac{f\!\left(\zeta'\right)}{1+V\!\left(\zeta'\right)}\cdot\frac{1+V\!\left(\zeta'\right)}{1+V(\zeta)}\right|P^{(2)}_t\!\left(\zeta,\mathrm{d}\zeta'\right)\notag\\[0.2em]
            &\leq \|f\|_V\frac{1+P_t^{(2)}V(\zeta)}{1+V(\zeta)}\notag\\[0.2em]
            &\leq \|f\|_V\frac{1+e^{\Lambda_1t}V(\zeta)}{1+V(\zeta)}\leq e^{\Lambda_1t}\|f\|_V,\label{eq: drift condition for small time t}
        \end{align}where $P_t^{(2)}V\leq e^{\Lambda_1t}V$ is given by the inverse-Lipschitz estimate. Writing $t=nT_*+r$ with $0\leq r<T_*$, the semigroup property combined with \eqref{eq: drift condition for nT*} and \eqref{eq: drift condition for small time t} implies
        \begin{equation*}
            \left\|P_t^{(2)}f\right\|_V\leq e^{\Lambda_1T_*}\alpha^n\|f\|_V.
        \end{equation*}Define $\kappa:=-\frac{1}{T_*}\ln\alpha>0$ and absorb the remaining fixed factors into $C$, we obtain
        \begin{equation}
            \left\|P_t^{(2)}f\right\|_V\leq Ce^{-\kappa t}\|f\|_V,\qquad \forall t\geq 0\label{eq: drift condition for all t}
        \end{equation}
        
        \item \label{step 2: proof of Theorem 1.1}  We denote the standard orthogonal basis for $L^2\!\left(\mathbb{T}^2\right)$ by $\left\{e_m(x)=e^{2\pi im\cdot x}:m\in\mathbb{Z}^2\right\}$ and let
        \begin{equation*}
            e_m^{(2)}(Z,x,y):=e_m(x)\overline{e_m}(y).
        \end{equation*}For a fixed constant $\beta:=\frac{\kappa}{14+2q}$, an initial vortex configuration $Z$ and $m,m'\neq \mathbf{0}$, we define 
        \begin{equation*}
            N_{m,m'}(Z):=\max\!\left\{n\geq 0:\left|\int e_m(x)e_{m'}(\phi_n(x))\mathrm{d}x\right|>e^{-\beta n}\right\}.
        \end{equation*}Since $e_{m'}^{(2)}$ has zero Lebesgue mean, \eqref{eq: drift condition for all t} gives
        \begin{align*}
            \mathbb{E}\left|\int e_m(x)e_{m'}(\phi_n(x))\mathrm{d}x\right|^2&=\left|\int e_m^{(2)}P_n^{(2)}e_{m'}^{(2)}\mathrm{d}x\mathrm{d}y\right|\\[0.2em]
            &\leq \int(1+V)\left\|P_n^{(2)}e_{m'}^{(2)}\right\|_V\mathrm{d}x\mathrm{d}y\\[0.2em]
            &\leq Ce^{-\kappa n}\|e_{m'}^{(2)}\|_V\!\int(1+V)\mathrm{d}x\mathrm{d}y\leq Ce^{-\kappa n},
        \end{align*}where we exploited the invariance of the Lebesgue measure under $P_t^{(2)}$ and $|e_m^{(2)}|\leq 1$. Chebyshev's inequality therefore yields the uniform tail estimate
        \begin{align}
            \mathbb{P}\!\left(N_{m,m'}(Z)>l\right)&\leq \sum_{n>l}\mathbb{P}\!\left(\left|\int e_m(x)e_{m'}(\phi_n(x))\mathrm{d}x\right|>e^{-\beta n}\right)\notag\\[0.2em]
            &\leq \sum_{n>l}e^{2\beta n}\mathbb{E}\left|\int e_m(x)e_{m'}(\phi_n(x))\mathrm{d}x\right|^2\notag\\[0.2em]
            &\leq\sum_{n>l}Ce^{(-\kappa+2\beta)n}\leq Ce^{(-\kappa+2\beta)l}\label{eq: tail estimate of N_m,m'}
        \end{align}Therefore, $N_{m,m'}$ is almost surely finite and 
        \begin{equation}
            \left|\int e_m(x)e_{m'}(\phi_n(x))\mathrm{d}x\right|\leq e^{\beta N_{m,m'}-\beta n}.\label{eq: correlation decay for e_m}
        \end{equation}

        Next, we introduce
        \begin{align*}
            M_0(Z)&:=\max\!\left\{|m|\vee |m'|:e^{\beta N_{m,m'}}>|m|\,|m'|\right\}\\
            \hat{D}(Z)&:=\max_{|m|,|m'|\leq M_0}e^{\beta N_{m,m'}},
        \end{align*}such that the estimates
        \begin{align}
            \mathbb{P}(M_0>M)&\leq 2\!\sum_{m,m':|m|>M}\mathbb{P}\!\left(e^{\beta N_{m,m'}}>|m||m'|\right)\leq CM^{-2(5+q)},\notag\\
            e^{\beta N_{m,m'}}&\leq \hat{D}|m||m'|\label{eq: upper bound for N_m,m'}
        \end{align}hold uniformly in $Z$. Let $f,g\in H^s$ be mean-zero observables with the Fourier expansions
        \begin{equation*}
            f=\sum_{m\in\mathbb{Z}^2\setminus\{\mathbf{0}\}}\hat{f}_me_m,\qquad g=\sum_{m\in\mathbb{Z}^2\setminus\{\mathbf{0}\}}\hat{g}_me_m.
        \end{equation*}
        Combining \eqref{eq: correlation decay for e_m} with \eqref{eq: upper bound for N_m,m'}, it follows that
        \begin{align}
            \left|\int f(x)g(\phi_n(x))\mathrm{d}x\right|&\leq e^{-\beta n}\sum_{m,m'}|\hat{f}_m||\hat{g}_{m'}|e^{\beta N_{m,m'}}\notag\\[0.2em]
            &\leq \hat{D}e^{-\beta n}\!\left(\sum_{m}|m||\hat{f}_m|\right)\!\left(\sum_{m'}|m'||\hat{g}_{m'}|\right)\notag\\[0.2em]
            &\leq \hat{D}e^{-\beta n}\left\|f\right\|_{H^3}\left\|g\right\|_{H^{3}}.\notag
        \end{align}To control the right-hand side in terms of $H^s$ norms, we invoke the following approximation lemma.
        \begin{lemma}[{\hspace{-0.05em}\cite[Lemma 7.1]{BBPS-mixing}}]\label{lem: control H^s' via H^s}
            Let $0<s<s'$ and let $h\in H^s$ be a mean-zero function. Then, for any $\epsilon>0$, there exists a mean-zero function $h_\epsilon\in H^{s'}$ such that the following estimates hold
            \begin{enumerate}[label=(\arabic*), ref=(\arabic*)]
                \item $\|h_\epsilon\|_{L^2}\lesssim\|h\|_{L^2}$.\vspace{0.2em}
                \item $\|h_\epsilon-h\|_{L^2}\lesssim \epsilon \|h\|_{H^s}$.\vspace{0.2em}
                \item $\|h_\epsilon\|_{H^{s'}}\lesssim_{s,s'}\epsilon^{-\frac{s'-s}{s}}\|h\|_{H^s}$.
            \end{enumerate}
        \end{lemma}
        Applying Lemma \ref{lem: control H^s' via H^s} with $s'=3$, we obtain
        \begin{align*}
            \left|\int f(x)g(\phi_n(x))\mathrm{d}x\right|&\leq \left|\int f_\epsilon(x)g_\epsilon(\phi_n(x))\mathrm{d}x\right|+\|g\|_{L^2}\|f-f_\epsilon\|_{L^2}+\|f_\epsilon\|_{L^2}\|g-g_\epsilon\|_{L^2}\\[0.2em]
            &\leq \hat{D}\!\left(e^{-\beta n}\epsilon^{-\frac{6-2s}{s}}+\epsilon\right)\|f\|_{H^s}\|g\|_{H^s}.
        \end{align*}Optimizing over $\epsilon$ by choosing $\epsilon=e^{-\frac{\beta s}{6-s}n}$, we deduce the discrete-time mixing estimate
        \begin{equation}
            \left|\int f(x)g(\phi_n(x))\mathrm{d}x\right|\leq C\hat{D}e^{-\gamma_{s,q}n}\|f\|_{H^s}\|g\|_{H^s},\label{eq: correlation decay for f,g at discrete time n}
        \end{equation}where $\gamma_{s,q}:=\frac{\beta s}{6-s}$. To generalize this to continuous time $t\geq 0$, we decompose $t = n + \hat{t}$ with $\hat{t} \in [0,1)$. Applying \eqref{eq: correlation decay for f,g at discrete time n} to $g\circ\phi_{\hat{t}}^{Z_n,\theta^n\omega}$ gives
        \begin{equation*}
            \left|\int f(x)g(\phi_t(x))\mathrm{d}x\right|\leq C\hat{D}e^{-\gamma_{s,q}n}\|f\|_{H^s}\|g\circ\phi_{\hat{t}}\|_{H^s}.
        \end{equation*}Since the velocity field $u$ admits uniform derivative bounds, there exists a uniform constant $C$ such that
        \begin{equation*}
            \sup_{\hat{t}\in[0,1)}\|g\circ\phi_{\hat{t}}\|_{H^s}\leq C\|g\|_{H^s}.
        \end{equation*}Hence, we have
        \begin{equation*}
            \left|\int f(x)g(\phi_t(x))\mathrm{d}x\right|\leq D(\omega,Z)e^{-\gamma_{s,q}t}\|f\|_{H^s}\|g\|_{H^s},
        \end{equation*}where $D(\omega,Z):=C\hat{D}(\omega,Z)e^{\gamma_{s,q}}$.

        Finally, it remains to estimate the moments of $D(\omega,Z)$ by analyzing $\hat{D}$. Although the definition of $\hat{D}$ depends on the initial configuration $Z$, its tail estimates are uniform in $Z$. We conclude that
        \begin{align*}
            \mathbb{E}[\hat{D}^q]&=\sum_{M=1}^\infty\mathbb{E}\!\left(\mathbf{1}_{\{M_0=M\}}\max_{m,m'\leq M}e^{q\beta N_{m,m'}}\right)\\[0.2em]
            &\leq \sum_{M=1}^\infty \!\left(\mathbb{P}(M_0=M)\right)^{1/2}\cdot\left\|\max_{m,m'\leq M}e^{q\beta N_{m,m'}}\right\|_{L^2}\\[0.2em]
            &\leq \sum_{M=1}^\infty M^{-(5+q)}\cdot\sum_{m,m'\leq M}\left\|e^{q\beta N_{m,m'}}\right\|_{L^2}\leq C\sum_{M=1}^\infty M^{-(1+q)},
        \end{align*}where the last inequality follows since $\left\|e^{q\beta N_{m,m'}}\right\|_{L^2}$ is uniformly bounded according to \eqref{eq: tail estimate of N_m,m'}.
    \end{enumerate}
\end{proof}

\vspace{0.5em}
\begin{appendix}
\section{Proof of Lemma \ref{lem: nondegeneracy lemma}}\label{sec: Proof of Lemma 2.13}
We complete the proof of Lemma \ref{lem: nondegeneracy lemma} by contradiction. For brevity, we only establish \eqref{eq:span one-point}-\eqref{eq:span Jocabian}, as \eqref{eq:span two-point} follows by an analogous argument.
\begin{proof}
    We begin with a common observation: if a real-analytic function $h\colon\mathbb T^2\to\mathbb R$ has zero spatial mean and satisfies
    \begin{equation*}
        \partial^{\alpha}h(r)=0,\qquad \forall\, |\alpha|\geq 1,
    \end{equation*}then $h\equiv 0$.

    To prove \eqref{eq:span one-point}, assume there exists $l\in\mathbb{R}^2$ such that
    \begin{equation*}
        l\cdot\partial^{\alpha}K_\epsilon(r)=0,\qquad\forall\, |\alpha|\geq 1.
    \end{equation*}Applying the preceding observation to $h(x):=l\cdot K_\epsilon(x)$ then yields $h\equiv 0$. By evaluating the corresponding Fourier coefficients, we obtain
    \begin{equation*}
        l\cdot k^\perp=0,\qquad \forall k\in\mathbb{Z}^2\setminus\{\mathbf{0}\}.
    \end{equation*}Choosing $k=(1,0)$ and $(0,1)$ confirms that $l=\mathbf{0}$.
    
    For the projective statement \eqref{eq:span projective}, suppose $l\in\mathbb{R}^2$ and $\beta=cv^\perp\in T_v\mathbb{S}^1$ satisfy
    \begin{equation*}
        l\cdot \partial^\alpha K_\epsilon(r)+\beta\cdot \Pi_v\partial^\alpha DK_\epsilon(r)v=0, \qquad \forall |\alpha|\geq 1.
    \end{equation*}Noting that $\beta\cdot\Pi_v=\beta\,\cdot$, the function 
    \begin{equation*}
        h(x):=l\cdot K_\epsilon(x)+\beta\cdot DK_\epsilon(x)v,
    \end{equation*}satisfies same analyticity and zero-mean conditions, and hence $h\equiv 0$. Comparison of Fourier coefficients gives, for every $k \neq \mathbf{0}$,
    \begin{equation*}
        i\,l\cdot k^\perp-2\pi(\beta\cdot k^\perp)(k\cdot v)=0.
    \end{equation*}Taking real and imaginary parts, we obtain
    \begin{equation*}
        l\cdot k^\perp=0,\qquad (\beta\cdot k^\perp)(k\cdot v)=c(k\cdot v)^2=0,\qquad \forall k\neq \mathbf{0},
    \end{equation*}which implies $l=\mathbf{0}$ and $c=0$.

    Finally, we prove \eqref{eq:span Jocabian}. The divergence-free property of $K_\epsilon$ gives
    \begin{equation*}
        \partial^\alpha DK_\epsilon(r)\in\mathfrak{sl}_2(\mathbb{R}),\qquad \forall \,|\alpha|\geq 1.
    \end{equation*}Suppose that $l\in\mathbb R^2$ and $C\in\mathfrak{sl}_2(\mathbb R)$ annihilate all pairs
    \begin{equation*}
        \left(\partial^\alpha K_\epsilon(r),\partial^\alpha DK_\epsilon(r)\right),\qquad |\alpha|\geq1,
    \end{equation*}we analogously deduce that
    \begin{equation*}
         h(x):=l\cdot K_\epsilon(x)+\operatorname{tr}\!\left(C^TDK_\epsilon(x)\right)\equiv 0.
    \end{equation*}Therefore, the Fourier expansion yields
    \begin{equation*}
        i\,l\cdot k^\perp-2\pi (k^\perp)^TCk=0,\qquad \forall k\neq\mathbf{0}.
    \end{equation*}The imaginary part gives $l=0$, while the real part implies
    \begin{equation*}
        (k^\perp)^TCk=0,\qquad \forall k\neq\mathbf{0}.
    \end{equation*}Writing 
    \begin{equation*}
        C=\begin{pmatrix}
            a & b\\
            c & -a
        \end{pmatrix},
    \end{equation*}and choosing $k=(1,0),(0,1)$, and $(1,1)$, we obtain $C=O_{2\times 2}$ as desired.
\end{proof} 
\end{appendix}

\noindent {\bf Acknowledgments.}
Yankai Shi is supported by NSFC 12231002. The author would like to express sincere gratitude to Prof. Yong Liu, Ziyu Liu, and Bin Tang for stimulating discussions, insightful comments, and continuous encouragement.

\normalem
\bibliographystyle{abbrv}
\bibliography{References}

\end{document}